\documentclass[%
onecolumn
, BCOR=1cm 
   , a4paper
]{mpi2015-cscpreprint}
\input{packages.tex}

\newcommand{\Gammacool}{\ensuremath{\Gamma_{\text{cool}}}}
\newcommand{\GammaN}{\ensuremath{\Gamma_\text{N}}}
\newcommand{\Gammaint}{\ensuremath{\Gamma_{\text{int}}}}

\newcommand{\Gammain}{\ensuremath{\Gamma_{\control}}}
\newcommand{\domain}{\ensuremath{\Omega}}
\newcommand{\domainl}{\ensuremath{\domain_l}}
\newcommand{\domains}{\ensuremath{\domain_s}}
\newcommand{\sysstatexstar}{\ensuremath{\sysstate^*}}
\newcommand{\Gammaintdiff}{\ensuremath{\Gamma_{\text{int},\Delta}(\ttime,\sysstatexstar)}}

\newcommand{\Gammainttx}{\ensuremath{{\Gamma}_{\text{int}}(\ttime,\sysstatexstar)}}
\newcommand{\GammaintRef}{\ensuremath{{\Gamma}_{\text{int,ref}}(\ttime, \sysstatexstar)}}

\newcommand{\Gammainx}[1][]{%
\ifthenelse{\isempty{#1}}{\Gammain}{\ensuremath{\Gamma_{\control,#1}}}}%
\newcommand{\Gammacoolx}[1][]{%
\ifthenelse{\isempty{#1}}{\Gammacool}{\ensuremath{\Gamma_{\text{cool},#1}}}}%
\newcommand{\Gammaoutx}[1][]{%
\ifthenelse{\isempty{#1}}{\ensuremath{\Gamma_{\cC}}}{\ensuremath{\Gamma_{\cC,#1}}}}%
\newcommand{\Tcool}{\ensuremath{\temp_\text{cool}}}

\newcommand{\kk}{\ensuremath{k}}
\newcommand{\ks}{\ensuremath{\kk_s}}
\newcommand{\kl}{\ensuremath{\kk_l}}
\newcommand{\Latent}{\ensuremath{\ell}}
\newcommand{\alphac}{\ensuremath{\alpha}}

\newcommand{\weight}{\ensuremath{\lambda}}

\providecommand{\temp}{} 
\renewcommand{\temp}{\ensuremath{\varTheta}}

\newcommand{\Vint}{\ensuremath{\V_{\text{int}}}}
\newcommand{\V}{\ensuremath{\Upsilon}}

\newcommand{\control}{\ensuremath{u}}
\newcommand{\controlk}{\ensuremath{\control\kiterind}}
\newcommand{\controlkprev}{\ensuremath{\control\kiterprevind}}
\newcommand{\controlkt}{\ensuremath{\tilde{\control}\kiterind}}

\newcommand{\controlref}{\ensuremath{\control_\text{ref}}}

\newcommand{\controlkref}{\ensuremath{\control_{\kiter,\text{ref}}}}
\newcommand{\controlf}{\ensuremath{\control_\cK}}
\newcommand{\intnormal}{\ensuremath{\normal_{\text{int}}}}
\newcommand{\normal}{\ensuremath{\boldsymbol{n}}}

\newcommand{\pert}{\ensuremath{\varphi}}
\newcommand{\controlx}[1][]{%
\ifthenelse{\isempty{#1}}{\control}{\ensuremath{\control_{#1}}}}%
\newcommand{\controltilde}[1][]{%
\ifthenelse{\isempty{#1}}{\control}{\ensuremath{\tilde{\control}_{#1}}}}%
\newcommand{\pertx}[1][]{%
\ifthenelse{\isempty{#1}}{\pert}{\ensuremath{\pert_{#1}}}}%
\newcommand{\genfunc}{\ensuremath{f}}

\newcommand{\genfuncFEM}{\ensuremath{f\hFEM}}

\newcommand{\integrand}{\ensuremath{\mathcal{Q}}}

\newcommand{\cM}{\ensuremath{\mathcal{M}}}
\newcommand{\cA}{\ensuremath{\mathcal{A}}}
\newcommand{\cB}{\ensuremath{\mathcal{B}}}
\newcommand{\cC}{\ensuremath{\mathcal{C}}}
\newcommand{\cK}{\ensuremath{\mathcal{K}}}
\newcommand{\cS}{\ensuremath{\mathcal{S}}}
\newcommand{\cI}{\ensuremath{\mathcal{I}}}
\newcommand{\cT}{\ensuremath{\mathcal{T}}}
\newcommand{\cO}{\ensuremath{\mathcal{O}}}
\newcommand{\bX}{\ensuremath{\mathbf{X}}}

\newcommand{\cBorig}{\ensuremath{\hat{\cB}}}

\newcommand{\cost}{\ensuremath{\mathcal{J}}}

\newcommand{\sysstate}{\ensuremath{x}}
\newcommand{\sysstatek}{\ensuremath{\sysstate\kiterind}}
\newcommand{\sysstatekt}{\ensuremath{\tilde{\sysstate}\kiterind}}

\newcommand{\sysstateref}{\ensuremath{\sysstate_\text{ref}}}
\newcommand{\sysstatekref}{\ensuremath{\sysstate_{\kiter,\text{ref}}}}

\newcommand{\sysstateDiff}{\ensuremath{\sysstate_\Delta}}

\newcommand{\sysstateDiffdot}{\ensuremath{\dot{\sysstate}_\Delta}}
\newcommand{\genstate}{\ensuremath{\zeta}}
\newcommand{\genstatek}{\ensuremath{\genstate\kiterind}}

\renewcommand{\output}{\ensuremath{y}}
\newcommand{\outputDiff}{\ensuremath{\output_\Delta}}

\newcommand{\nstate}{\ensuremath{n}}
\newcommand{\ninput}{\ensuremath{m}}
\newcommand{\noutput}{\ensuremath{p}}
\newcommand{\nerror}{\ensuremath{n_e}}
\newcommand{\nin}{\ninput}
\newcommand{\nout}{\noutput}

\newcommand{\nLt}{\ensuremath{s\time}}
\newcommand{\nLtzero}{\ensuremath{s(\Starttime)}}
\newcommand{\ord}{\ensuremath{\wp}}
\newcommand{\nord}{\ensuremath{n_{\ord}}}

\newcommand{\ntimestep}{\ensuremath{n_t}}
\newcommand{\kiter}{\ensuremath{k}}
\newcommand{\kiterind}{\ensuremath{_\kiter}}
\newcommand{\kiterprevind}{\ensuremath{_{\kiter - 1}}}
\newcommand{\kiternextind}{\ensuremath{_{\kiter + 1}}}
\newcommand{\kiterindbwd}{\ensuremath{_{\ntimestep-\kiter}}}
\newcommand{\kiterindT}{\ensuremath{^{\phantom{\transpose}}\kiterind}}
\newcommand{\kiterindTbwd}{\ensuremath{^{\phantom{\transpose}}\kiterindbwd}}
\renewcommand{\dim}{\ensuremath{d}}

\newcommand{\ttime}{\ensuremath{t}}
\renewcommand{\time}{\ensuremath{(\ttime)}}
\newcommand{\Endtime}{\ensuremath{\ttime_\text{end}}}
\newcommand{\Starttime}{\ensuremath{\ttime_0}}
\newcommand{\tk}{\ensuremath{\ttime\kiterind}}

\newcommand{\tkprev}{\ensuremath{\ttime\kiterprevind}}

\newcommand{\tknextref}{\ensuremath{\ttime_{\kiter + 1, \text{ref}}}}

\newcommand{\timek}{\ensuremath{(\tk)}}

\newcommand{\timeinterval}{\ensuremath{[\Starttime,\Endtime]}}
\newcommand{\timeintervalk}{\ensuremath{[\tkprev,\tk]}}
\newcommand{\thetak}{\ensuremath{\Theta}}
\newcommand{\diff}[1]{\ensuremath{\mathrm{d}#1}}
\newcommand{\timediff}{\ensuremath{\frac{\mathrm{d}}{\mathrm{d}t}}}
\newcommand{\texp}[1]{\ensuremath{\operatorname{exp}\left(#1\right)}}
\newcommand{\mess}{\textsf{M.E.S.S.}}
\newcommand{\mmess}{\mbox{\textsf{M-}\mess{}}}
\newcommand{\transpose}{\ensuremath{^{\textsf{T}}}}
\newcommand{\invtranspose}{\ensuremath{^{-\textsf{T}}}}
\newcommand{\discreteset}{\ensuremath{\mathscr{T}}}
\newcommand{\timesteps}{\discreteset}

\newcommand{\timestepsfwd}{\ensuremath{\discreteset_{\text{fwd}}}}
\newcommand{\timestepsfwdref}{\ensuremath{\discreteset_{\text{fwd}}^{\text{ref}}}}
\newcommand{\meshcells}{\ensuremath{\mathscr{Q}{\hFEM}}}
\newcommand{\meshcell}{\ensuremath{Q}}
\newcommand{\timestep}{\ensuremath{\tau}}
\newcommand{\timestepk}{\ensuremath{\timestep\kiterind}}

\newcommand{\timestepknext}{\ensuremath{\timestep\kiternextind}}
\newcommand{\hFEM}{\ensuremath{^h}}
\newcommand{\FEMspace}{\ensuremath{\mathscr{V}}}
\newcommand{\sysstatespace}{\ensuremath{\mathscr{V}_\sysstate}}
\newcommand{\inspace}{\ensuremath{\mathscr{V}_{\control}}}

\newcommand{\outspace}{\ensuremath{\mathscr{V}_{\text{out}}}}
\newcommand{\timespace}{\ensuremath{\mathscr{I}}}

\newcommand{\indicatork}{\ensuremath{\mathcal{I}\kiterind}}
\newcommand{\indicatorke}{\ensuremath{\indicatork^e}}
\newcommand{\indicatorku}{\ensuremath{\indicatork^{\Delta\control}}}
\newcommand{\indicatorkdtu}{\ensuremath{\indicatork^{\dot\control}}}
\newcommand{\indicatorkt}{\ensuremath{\delta}}
\newcommand{\indicatorfunc}{\ensuremath{\mathcal{I}(\tk, \genstatek, \controlk)}}
\newcommand{\indicatorup}{\ensuremath{\overline{\indicatorkt}}}
\newcommand{\indicatordown}{\ensuremath{\underline{\indicatorkt}}}
\newcommand{\retrystep}{\ensuremath{\textsf{RetryStep}}}
\newcommand{\tol}{\ensuremath{\textsf{TOL}}}
\newcommand{\rstep}{\ensuremath{r}}
\newcommand{\gammastep}{\ensuremath{\gamma}}
\newcommand{\pickfkt}{\ensuremath{\chi}}
\newcommand{\errDRE}{\ensuremath{e_{\text{DRE}}}}

\newcommand{\avg}{\ensuremath{\mathop{\mathrm{avg}}}}
\newlength\figureheight
\newlength\figurewidth
\newcounter{mymac@matlab}
  \newcommand{\matlab}{MATLAB%
   \ifnum\value{mymac@matlab}<1%
   \textsuperscript{\textregistered}%
   \setcounter{mymac@matlab}{1}%
   \fi%
  }
\newcommand{\fenics}{\textsf{FEniCS}}
\newcommand{\numpy}{\textsf{NumPy}}
\newcommand{\scipy}{\textsf{SciPy}}
\newenvironment{blockmatrix}{%
  \left[\vphantom{n}%
  \vcenter\bgroup\hbox\bgroup
  \tikzpicture[
    x=1\baselineskip,
    y=1\baselineskip,
    remember picture,
  ]%
}{%
  \endtikzpicture
  \egroup
  \egroup
  \vphantom{n}\right]%
} 

\newcommand*{\mblock}[1][white]{%
  \mblockaux{#1}%
}
\def\mblockaux#1(#2,#3)#4(#5,#6){%
  \draw[fill={#1},draw=white,line width=1pt]
  let \p1=(#2,#3),
      \p2=(#5,#6),
      \p3=(#2+#5,#3+#6),
      \p4=(#2+#5/2,#3+#6/2)
  in
    (\p1) rectangle (\p3)
    (\p4) node {$#4$}
  ;%
}

\def\lblockaux#1(#2,#3)#4(#5,#6)(#7){%
  \draw[fill={#1},draw=white,line width=1pt]
  let \p1=(#2,#3),
      \p2=(#5,#6),
      \p3=(#2+#5,#3+#6),
      \p4=(#2+#5/2,#3+#6/2)
  in
    (\p1) rectangle (\p3)
    (\p4) node(#7) {$#4$}
  ;%
}
\newcommand\restr[2]{{
  \left.\kern-\nulldelimiterspace 
  #1 
  \vphantom{\big|} 
  \right|_{#2} 
  }}
\makeatletter 
\def\ifEmpty#1{\def\@temp{#1}\ifx\@temp\@empty} 
\makeatother
\newcommand{\nrm}[2][]{\ensuremath{\left\lVert #2 \right\rVert\ifEmpty{#1}\else_{#1}\fi}}

\definecolor{uncontroled}{RGB}{217,95,2}
\definecolor{desiredpre}{RGB}{102,166,30}
\colorlet{desired}{desiredpre!20!black}
\definecolor{controled}{RGB}{117,112,179}
\definecolor{extracolor1}{RGB}{27,158,119}
\definecolor{extracolor2pre}{RGB}{230,171,2}
\colorlet{extracolor2}{extracolor2pre!95!black}
\definecolor{extracolor3}{RGB}{231,41,138}
\colorlet{inputcolor}{orange}
\definecolor{outputcolor}{HTML}{33a5c3}
\colorlet{pertcolor}{red!70!black}
\colorlet{bdf1color}{uncontroled}
\colorlet{bdf2color}{desired}
\colorlet{bdf3color}{controled}
\colorlet{bdf4color}{extracolor1}
\colorlet{split1color}{extracolor2}
\colorlet{split2color}{extracolor3}
\colorlet{nofeedbackcolor}{pertcolor}
\colorlet{fterr8color}{blue}
\colorlet{fterr9color}{green}
\colorlet{fterr10color}{orange}
\colorlet{ftu2color}{red}
\colorlet{ftdtu2color}{cyan}
\tikzset{cross/.style={cross out, draw, 
         minimum size=2*(#1-\pgflinewidth), 
         inner sep=0pt, outer sep=0pt}}

\crefname{line}{Line}{Lines}
\crefalias{AlgoLine}{line}%
\Crefname{algocf}{Algorithm}{Algorithms}
\makeatletter
\let\cref@old@stepcounter\stepcounter
\def\stepcounter#1{%
  \cref@old@stepcounter{#1}%
  \cref@constructprefix{#1}{\cref@result}%
  \@ifundefined{cref@#1@alias}%
    {\def\@tempa{#1}}%
    {\def\@tempa{\csname cref@#1@alias\endcsname}}%
  \protected@edef\cref@currentlabel{%
    [\@tempa][\arabic{#1}][\cref@result]%
    \csname p@#1\endcsname\csname the#1\endcsname}}
\makeatother

\begin{document}

\title{Numerical integrators for non-autonomous closed-loop systems with
  Riccati-feedback}

\author[$\ast$]{Bj\"orn Baran}
\affil[$\ast$]{Max Planck Institute for Dynamics of Complex Technical Systems, Sandtorstr. 1, 39106 Magdeburg.\authorcr%
  \email{baran@mpi-magdeburg.mpg.de}, \orcid{0000-0001-6570-3653}
}
\author[$\ast\ast$]{Peter Benner}
\affil[$\ast\ast$]{Max Planck Institute for Dynamics of Complex Technical Systems, Sandtorstr. 1, 39106 Magdeburg.\authorcr%
  \email{benner@mpi-magdeburg.mpg.de}, \orcid{0000-0003-3362-4103}
}
\author[$\dagger$]{Jens Saak}
\affil[$\dagger$]{Max Planck Institute for Dynamics of Complex Technical Systems, Sandtorstr. 1, 39106 Magdeburg.\authorcr%
  \email{saak@mpi-magdeburg.mpg.de}, \orcid{0000-0001-5567-9637}
}
\author[$\ddagger$]{Tony Stillfjord}
\affil[$\ddagger$]{Lund University, Box 117, SE-221 00, Lund, Sweden.\authorcr%
  \email{tony.stillfjord@math.lth.se}, \orcid{0000-0001-6123-4271}
}

\shortdate{}
  



%
\shorttitle{Num. integr. for non-auton. closed-loop syst. with Riccati-feedback}
\shortauthor{B. Baran, P. Benner, J. Saak, T. Stillfjord}

\abstract{%
By computing a feedback control via the linear quadratic regulator (LQR) approach
and simulating a non-linear non-autonomous closed-loop system using this feedback, we
combine two numerically challenging tasks.
For the first task, the computation of the feedback control, we use the
non-autonomous generalized differential Riccati equation (DRE), whose
solution determines the time-varying feedback gain matrix.
Regarding the second task, we want to be able to simulate non-linear closed-loop systems for
which it is known that the regulator is only valid for sufficiently
small perturbations. Thus, one easily runs into numerical issues
in the integrators when the closed-loop control varies greatly. 
For these systems, e.g., the A-stable implicit Euler methods fails.\newline
On the one hand, we implement non-autonomous versions 
of splitting schemes and BDF methods for the solution of our 
non-autonomous DREs.
These are well-established DRE solvers in the autonomous case.  
On the other hand, to tackle the numerical issues in the
simulation of the non-linear closed-loop system, we apply a
fractional-step-theta scheme with time-adaptivity tuned specifically
to this kind of challenge. That is, we additionally base the time-adaptivity on the
activity of the control. We compare this approach to the more
classical error-based time-adaptivity.\newline
We describe techniques to make these two tasks computable in a reasonable amount of time
and are able to simulate closed-loop systems
with strongly varying controls, while avoiding numerical issues.
Our time-adaptivity approach requires fewer time steps than the error-based alternative and is more reliable.
}

\novelty{%
The novelty of this work is twofold:
\begin{itemize}
\item Until now, numerical methods for non-autonomous differential Riccati
      equations have been restricted to very specific assumptions on the
      time-dependent matrices. We propose new low-rank splitting schemes and
      backward differentiation formulas to numerically solve general
      non-autonomous differential Riccati equations.
\item We overcome numerical issues that occur during the simulation of
      the corresponding closed-loop system by utilizing the
      fractional-step-theta algorithm combined with a time-adaptive strategy
      that is specifically adapted to quickly varying controls.
\end{itemize}
}

\maketitle

\section{Introduction}\label{sec:into}

The combination of the two tasks, to compute a feedback control for a (possibly) non-linear non-autonomous system
and to apply this control in a simulation of the closed-loop system, consists of three steps.
First, we simulate the (open-loop) non-linear problem, (possibly) with a
feedforward-control, to assemble the reference solution trajectory and data for
the feedback computation.
Then, we linearize the equations around the reference trajectory and
apply the linear-quadratic regulator (LQR) approach, to the linearized system, to receive the feedback gain.
In the third step, we apply the feedback gain to the non-linear problem during a closed-loop simulation.
The novelty of this approach lies in the application of this framework to
large-scale non-autonomous systems.
This can include problems with free boundaries and moving interfaces.

The second step, the computation of a feedback control with the LQR approach 
on a finite time-horizon, results in the necessity
to solve a differential Riccati equation (DRE)~\cite{Rei72, BitG91}.
In fact, DREs arise for other applications as well
and have been studied in great detail~\cite{AboFIetal03}.
We focus in particular on the generalized DREs~\cite{KunM90} and DREs 
resulting from differential-algebraic equations (DAEs)~\cite{KunM89}.
The algebraic constraint equations can arise, e.g., from divergence-free conditions in convection-diffusion equations or from coupling constraints 
with free boundaries and moving interfaces.

In order to numerically solve the corresponding DRE
as well as the non-linear problem in step one and three,
we first semi-discretize in space.
Thereby, the DRE becomes matrix-valued with large-scale, but sparse, coefficients and 
a solution that, in most cases, can be approximated sufficiently accurate 
by a low-rank factorization~\cite{Sti18}.
There is a range of methods available to numerically solve this kind of DREs.
To enumerate several, splitting schemes 
can be applied to DREs~\cite{Sti15,Sti18a,Sti15a,OstPW19,MenOPetal18}.
Further, there are the backward differentiation formulas (BDF)~\cite{BenM04, Men07, LanMS15, BenM18} as well as Rosenbrock and Peer methods~\cite{LanMS15,Men07,Lan17,BenL18}.
Recently, also Krylov subspace methods have been applied to DREs~\cite{KirSi19,GueHJetal18,BehBH21,KosM20}
as well as exponential integrators~\cite{LiZL20}
and an all-at-once space-time approach~\cite{BreDS21} (computing a low-rank tensor format solution).

The just mentioned methods are mainly focused on autonomous DREs.
In contrast to this, DREs are especially challenging when the coefficient matrices are time-dependent. 
In detail, these non-autonomous DREs are particularly challenging for problems with free boundaries~\cite{BaePS13, BarBHetal18a} and moving interfaces like 
fluid structure interaction (FSI)~\cite{Fai17} and the Stefan problem~\cite{BarBHetal18a,BarBS22a}.
Addressing such problems, the numerical solution of large-scale non-autonomous DREs has not been studied in a lot of detail so far.
For instance, in~\cite{Men07, BenM18}, the authors formulate non-autonomous BDF
and Rosenbrock methods without going into detail about the implementation or performance.
Moreover, non-autonomous BDF and Rosenbrock methods are used in~\cite{LanMS15,Lan17} as well as non-autonomous Peer methods in~\cite{BenL18}, but are restricted 
to the case of constant matrices scaled by time-dependent scalar functions. 
The applications we consider have no restriction on the non-autonomy of the matrices.
Consequently, in this manuscript, we extend the existing BDF methods (and to a limited extent also splitting methods) for solving large-scale fully non-autonomous DREs.
Further, we give details on how to correctly incorporate the weight parameter from the cost functional into the DRE and provide a numerical implementation in the software package \mmess~\cite{SaaKB-mmess-all-versions}.

In addition to these contributions, this manuscript focuses on
solving closed-loop problems with free boundaries and moving interfaces.
This is a very broad and active field of research,
see e.g.,~\cite{KogDK19,LuF20,JinWN20,CvoRD20,KogMCetal20},
and there are numerous numerical methods available.
During the simulation of problems with moving interfaces, 
we can observe numerical issues, which are similar to the ones in~\cite{FaiW18}.
This is especially true for closed-loop problems with strongly varying controls.
In this setting, common methods like the implicit Euler method and the trapezoidal rule
run the risk to break down as in~\cite{FaiW18}.
To overcome this, our method of choice is the fractional-step-theta algorithm~\cite{BriGP87} combined with adaptive time stepping.
It has successfully been applied for Navier-Stokes equations~\cite{MeiR14, MeiR15}
as well as for FSI problems~\cite{FaiW18,Wic11,RicW15}.
Part of this method can be classical error estimates for the time-adaptivity, e.g., heuristics, based on models of the actual error~\cite{GusLS88, Tur99, JohR10},
like the one that is used in~\cite{FaiW18}.
Others are residual-based like the dual weighted residual, which yields the best results in~\cite{FaiW18}.
Unfortunately, this method has high computational costs and, thus, is not feasible for our purpose.
Therefore, our approach is to modify an error-based heuristic 
and adapt it to the case of closed-loop systems with quickly varying controls.

For specific problems that have moving inner boundaries as well, 
there are time-adaptive strategies available that differ from ours.
One example is an adaptive time stepping applied for a two-phase flow
with a proportional–integral–derivative (PID) control in~\cite{AkaG17}.
Another example is~\cite{SonUS15}, where temporal and spatial discretizations 
are used that are adaptive with respect
to Courant-Friedrichs-Lewy and non-self-intersection conditions for a dendritic Stefan problem.

\paragraph{Structure of the Manuscript.}
In \cref{sec:problem}, we state a general definition of the (possibly)
non-linear non-autonomous closed-loop system we consider.
To design a feedback controller for this system, we choose the LQR approach and formulate the non-autonomous DRE in \cref{sec:LQR}.
To solve this DRE, the non-autonomous versions of the BDF methods and the splitting schemes are discussed in \cref{sec:DRE}.
With this solution, we can compute a feedback control, which is then used in the closed-loop system.
For a robust numerical solver of this system,
we describe several time stepping methods for closed-loop integrators, e.g., a time-adaptive fractional-step-theta algorithm
with several indicators in \cref{sec:time}.
We demonstrate and compare the performance of the presented methods and their behavior 
with respect to different control-parameter choices in \cref{sec:numeric}.
Additionally, we specify how to access the codes and data for our experiments in~\cite{LQR_Stefan_codes}
and summarize our conclusions in \cref{sec:conclusion}.

\paragraph{Notation.}
In most equations, we omit the time-dependence \time{}, the spatial dependence $(\sysstate)$, or the combination of both $(\ttime, \sysstate)$.
This is intended to improve the readability of the equations.

\subsection{Problem Definition}\label{sec:problem}
Our goal is to derive and apply a feedback control \controlf{} for possibly non-linear non-autonomous systems.
This section describes the corresponding control problem in an abstract form such that it can be applied to different dynamical systems.

We define this abstract, controlled, dynamical system
\begin{equation}\label{eq:system}
  \dot\sysstate = \genfunc(\sysstate, \control)
\end{equation}
on a finite time horizon $\ttime \in \timespace = \timeinterval$
and a domain $\domain \subset \mathbb{R}^\dim$ of dimension $\dim > 0$.
The domain has a fixed outer shape and can contain moving interfaces or free inner boundaries.
However, we expect that our methods work as well with a time-dependent domain that is finite and 
with a volume that is uniformly bounded away from zero.
Additionally, the function $\genfunc$ is possibly non-linear in \sysstate\time{}.
This state $\sysstate \in \sysstatespace$  as well as the control $\control \in \inspace$ are unknowns, 
with the function spaces \sysstatespace{} and \inspace{} being suitable spaces chosen with 
respect to the underlying dynamical problem.

Further, we assume that a reference pair $(\sysstateref, \controlref) \in \sysstatespace\times\inspace$ is given, which is a solution to \cref{eq:system}.
In particular, this represents the desired state trajectory of the system, which we want to stabilize with the feedback control \(\controlf\).
That means, in this scenario, that the resulting input is given as $\control = \controlref + \controlf$.
To clarify, the additional feedback control is necessary because
we assume that the system~\labelcref{eq:system} is influenced by model inaccuracies or perturbations,
which cause the state to deviate from the desired trajectory.
This deviation is expressed as 
\begin{equation*}
  \sysstateDiff = \sysstate - \sysstateref, \qquad 
  \controlf(\ttime, \sysstateDiff\time) = \control - \controlref. 
\end{equation*}
Consequently, we reformulate system~\labelcref{eq:system} in closed-loop form
\begin{equation}\label{eq:closedloop}
  \sysstateDiffdot = \genfunc(\sysstateDiff, \controlf).
\end{equation}
With this setup, the desired state of system~\labelcref{eq:closedloop} is the zero state.

In order to measure the deviation of the state in a performance index,
we denote the function that maps the state to the output~\outputDiff\time{} of the system by $\cC\colon\sysstatespace \to \outspace$,
where \outspace{} is a suitable space:
\begin{equation*}
  \outputDiff = \cC\sysstateDiff.
\end{equation*}

These output and control costs are measured in a performance index $\cost\colon \sysstatespace\times\inspace \to\mathbb{R}_{\ge0}$, defined by
\begin{equation}\label{eq:cost}
  \cost(\sysstateDiff,\controlf) = \int^{\Endtime}_{\Starttime} \nrm{\outputDiff}^2 + \weight\nrm{\controlf}^2 \diff{t} 
  + \sysstateDiff\transpose(\Endtime)S\sysstateDiff(\Endtime),
\end{equation}
where $S$ is positive semidefinite and $\weight{} > 0$ is a regularization parameter.
A smaller \weight{} makes the problem numerically more challenging and allows the feedback control greater activity.
As a result, minimizing the performance index corresponds to
the objective to stabilize the system~\labelcref{eq:closedloop} to zero.

To perform this task computationally,
we use numerical methods to solve the system~\labelcref{eq:closedloop} and
compute a feedback control $\ttime \mapsto \controlf(\ttime, \sysstateDiff\time)$,
thus approximating an optimal pair $(\sysstate^*, \control^*) \in \sysstatespace\times\inspace$ 
that minimizes the performance index $\cost(\sysstate^*, \control^*)$ 
with the LQR approach.

\subsection{Linear Quadratic Regulator}\label{sec:LQR}
The optimal feedback control for a linear control system in state-space formulation with a quadratic cost functional is given by the linear quadratic regulator, see e.g.~\cite{Loc01}.
Moreover, non-linear systems can be stabilized with this method if the deviation from the desired trajectory is sufficiently small~\cite[Section 8.5]{Son98}.
We use this approach because it is well studied for related types of problems, e.g., convection diffusion equations~\cite{Wei16}, and demonstrates promising performance for these.
For the non-linear case, first, we linearize system~\labelcref{eq:closedloop} around the reference pair $(\sysstateref, \controlref)$.
Next, we apply the LQR approach to the linearization.
Then, the derived feedback control can be applied to the non-linear system.

To begin with, the linearization of system~\labelcref{eq:closedloop} can always be written in the form
\begin{equation}\label{eq:sss}\begin{aligned}
  \cM\dot\sysstate\hFEM 	&= \cA\sysstate\hFEM + \cBorig\control\hFEM, \\
  \output\hFEM 			&= \cC\sysstate\hFEM.
\end{aligned}\end{equation}
In detail, the system is considered to be semi-discretized in space with the sparse square matrices 
$\cA\time, \cM\time \in \mathbb{R}^{\nstate \times \nstate}$, 
the input matrix $\cBorig\time \in \mathbb{R}^{\nstate \times \ninput}$, and
the output matrix $\cC\time \in \mathbb{R}^{\noutput \times \nstate}$.
These matrices are possibly time-dependent
and we assume \cM\time{} to be differentiable and uniformly non-singular.
Further, $\sysstate\hFEM \in \sysstatespace\hFEM \subset \timespace \times \mathbb{R}^\nstate$, 
$\control\hFEM \in \inspace\hFEM \subset \timespace \times \mathbb{R}^\nin$, and 
$\output\hFEM \in \outspace\hFEM \subset \timespace \times \mathbb{R}^\nout$ with $\sysstatespace\hFEM \subset \sysstatespace,
\inspace\hFEM \subset \inspace$, and $\outspace\hFEM \subset \outspace$ being semi-discrete subspaces.
Additionally, we assume that $\ninput, \noutput \ll \nstate$, 
such that low-rank methods are applicable in the next step.

The main task for computing an optimal feedback control is to solve 
a matrix-valued non-autonomous generalized differential Riccati equation (DRE), given by
\begin{equation}\label{eq:DREpre}
  \begin{aligned}
      - \timediff{(\cM\transpose \bX \cM)} &= \cC\transpose\cC+\cA\transpose\bX\cM+\cM\transpose \bX\cA
      -\cM\transpose\bX\cB\cB\transpose\bX\cM, \\
     \cM(\Endtime)\transpose \bX(\Endtime) \cM(\Endtime) &= S.
  \end{aligned}
\end{equation}
In contrast to this general non-autonomous problem, the DRE with constant mass matrix, i.e., \(\cM\time{}\equiv\hat\cM\), 
is more common in literature~\cite{Men07,BenM18,LanMS15,Lan17,BenL18} and
requires no special treatment of the time-derivative at the left-hand side of \cref{eq:DREpre}.
Here, this time-derivative can be reformulated using the chain rule:
\begin{equation*}
 	- \timediff{(\cM\transpose \bX \cM)} = - \dot\cM\transpose \bX \cM - \cM\transpose \dot\bX \cM - \cM\transpose \bX \dot\cM .
\end{equation*}
With a constant mass matrix, the terms on the left and right vanish and \cref{eq:DREpre} simplifies to
\begin{equation}\label{eq:DREconst}
  \begin{aligned}
      - \hat\cM\transpose \dot\bX \hat\cM &= \cC\transpose\cC+\cA\transpose\bX\hat\cM+\hat\cM\transpose \bX\cA
      -\hat\cM\transpose\bX\cB\cB\transpose\bX\hat\cM, \\
     \hat\cM(\Endtime)\transpose \bX(\Endtime) \hat\cM(\Endtime) &= S.
  \end{aligned}
\end{equation}
Otherwise, by subtracting the two terms containing $\dot\cM\time$ from \cref{eq:DREpre} we get a DRE 
that contains the time-derivative of \cM\time{}:
\begin{equation}\label{eq:DRE}
  \begin{aligned}
      - \cM\transpose\dot \bX \cM &= \cC\transpose\cC+(\dot\cM + \cA)\transpose\bX\cM+\cM\transpose \bX(\dot\cM + \cA)
      -\cM\transpose\bX\cB\cB\transpose\bX\cM, \\
     \cM(\Endtime)\transpose \bX(\Endtime) \cM(\Endtime) &= S.
  \end{aligned}
\end{equation}
Here, the coefficient matrices are time-dependent and so is the $\mathbb{R}^{\nstate \times \nstate}$-valued solution $\bX\time$.
Note that the input matrix $\cB\time = \frac{1}{\sqrt{\weight}}\cBorig\time$ is scaled with the weight factor from the performance index~\labelcref{eq:cost}. 

With the solution of \cref{eq:DRE}, we can compute the feedback gain matrix
\begin{equation*}
  \cK\ = \frac{1}{\weight}\cBorig\transpose\bX\cM = \frac{1}{\sqrt{\weight}}\cB\transpose\bX\cM.
\end{equation*}
Then, we obtain the feedback control $\controlf\hFEM(\ttime, \sysstateDiff\hFEM\time)$ that minimizes the performance index~\labelcref{eq:cost}
by applying $-\cK\time$ to the semi-discrete state deviation $\sysstateDiff\hFEM\time$ as in~\cite{Loc01,Son98,BenLP08}:
\begin{equation}\label{eq:feedbackcontrol}
  \controlf\hFEM = -\cK \sysstateDiff\hFEM.
\end{equation}
This semi-discrete feedback control can then be applied to the non-linear system~\labelcref{eq:closedloop}.

One main contribution of this manuscript are the numerical methods to efficiently solve the DRE~\labelcref{eq:DRE}
and, thus, compute the feedback gain \cK\time{}.
Such methods, which can also treat the time-dependent coefficients, are presented in the next section.

\section{DRE Solver}\label{sec:DRE}
In this section, we present methods to numerically compute a low-rank
approximation of the solution $\bX\time \in \mathbb{R}^{\nstate \times \nstate}$ of the DRE~\labelcref{eq:DRE}.
For the autonomous DRE~\labelcref{eq:DREconst}, i.e., when all matrix coefficients are constant in time, many efficient low-rank methods exist, as we have discussed in the introduction.
These methods are mostly tailored for the autonomous DRE and not yet suited to handle time-dependent coefficients.
Here, we introduce the BDF methods and splitting schemes adapted to the non-autonomous DRE~\labelcref{eq:DRE}.

To this end, we introduce an equidistant temporal discretization consisting of the $\ntimestep + 1$ grid points $\ttime\kiterind$ defined by
\begin{equation*}
  \timestep = \frac{\Endtime - \Starttime}{\ntimestep}, \qquad
  \timesteps = \{ \ttime\kiterind = \Starttime + \kiter\cdot\timestep, \kiter = 0,\ldots, \ntimestep \}.
\end{equation*}
The DRE is solved backwards in time starting at \Endtime.
For reasons of better readability and for consistency 
with previously published material on the autonomous case, 
we perform a change of variables ${\ttime \to \Endtime + \Starttime - \ttime}$ 
such that we can formulate the methods forward in time. 
The DRE to be solved is therefore
\begin{equation}
  \label{eq:DRE_fwd}
  \begin{aligned}
    \cM\transpose \dot{\bX} \cM
    &= \cC\transpose\cC
     + (\cA + \dot{\cM})\transpose \bX \cM
     + \cM\transpose \bX (\cA + \dot{\cM})
     - \cM\transpose\bX\cB\cB\transpose\bX\cM, \\
     \cM(\Starttime)\transpose \bX(\Starttime) \cM(\Starttime) &= \bX_0 = S.
  \end{aligned}
\end{equation}
The DRE solvers in this section then efficiently approximate $\bX\kiterindT \approx \bX\timek$,
representing $\bX\kiterindT$ by a low-rank factorization.

\subsection{Low-rank Methods}\label{sec:low-rank}
Motivated by~\cite{Sti18}, with $\ninput, \noutput \ll \nstate$, we assume the solution
of the DRE, \bX\time{}, to have a small (numerical) rank $\nLt \ll \nstate$.
Even if the matrices $\cA\time, \cM\time \in \mathbb{R}^{\nstate \times \nstate}$
are sparse, at each time \time, the solution $\bX\time \in \mathbb{R}^{\nstate \times \nstate}$ is a dense matrix. 
However, since it has a low numerical rank, it can be well approximated by
\begin{equation}\label{eq:LR}
  \bX \approx L D L\transpose,
\end{equation}
where the low-rank factors 
$L\time \in \mathbb{R}^{\nstate \times \nLt}$ and $D\time \in \mathbb{R}^{\nLt \times \nLt}$
have the rank $\nLt \ll \nstate$~\cite{LanMS15}.
With this low-rank factorization, just the low-rank factors 
$L\time$ and $D\time$
need to be stored. 
Consequently, the memory requirement for the storage of the solution reduces to $\cO(\nLt\nstate+\nLt^2)$
instead of $\cO(\nstate^2)$, per time step.
This factorization technique is used by the methods we introduce in this section.

\paragraph{Non-autonomous BDF Method}
We choose the BDF methods because the existing implementations can be adapted 
to the non-autonomous case relatively straight forward since 
the original method is, in theory, not restricted to autonomous data.
However, the existing BDF methods for matrix-valued DREs like~\cite[Algorithm 3.2]{LanMS15}
are tuned for constant data.
Nevertheless, we can extend these BDF methods to solve the non-autonomous DRE~\labelcref{eq:DRE} with a feasible amount of extra computational cost. 
This results in \cref{algo:BDF} formulated forward in time.
Note that the actual numerical implementation is backward in time.
An equivalent formulation of the non-autonomous BDF methods backward in time can be found in~\cite{BarBS22a}.
\begin{algorithm}[t]
	\caption{Non-autonomous low-rank factor BDF method of order \ord{}}%
	\label{algo:BDF}
	\hspace*{\algorithmicindent} \textbf{Input} $\cA\time,\cM\time,\dot\cM\time,\cB\time,\cC\time,\weight,\timesteps, \ord, L_0,\ldots,L_{\ord-1}, D_0,\ldots,D_{\ord-1}$ \\
	\hspace*{\algorithmicindent} \textbf{Output} $\cK\kiterind, \kiter = 0,\ldots,\ntimestep - 1$
\begin{algorithmic}[1]
	\State Invert the direction of time in all matrix functions, e.g.\ $\cA\time \to \cA(\Endtime  + \Starttime - \ttime)$\label{line:BDF0}
	\State $\cK\kiterindTbwd = \frac{1}{\sqrt{\weight}}\cB\timek\transpose L\kiterindT D\kiterindT L\kiterind\transpose\cM\timek, \quad \kiter = 1,\ldots,\ord - 1$\label{line:BDF01}
	\For{$\kiter = \ord,\ldots,\ntimestep$}\label{line:BDF1}
	  \State $\cA\kiterind = \timestep\beta(\dot\cM\timek + \cA\timek) - \frac{1}{2}\cM\timek$\label{line:BDF3}
	  \State $\cM\kiterind = \cM\timek$\label{line:BDF4}
	  \State $\cB\kiterind = \sqrt{\timestep\beta}\cB\timek$\label{line:BDF5}
	  \State $\cC\kiterind\transpose = \left[\cC\timek\transpose, \cM\kiterind\transpose L_{\kiter - 1}, \ldots,
	  \cM\kiterind\transpose L_{\kiter - \ord}\right]$\label{line:BDF6}
	  \State $\cS\kiterind =  
	    \begin{blockmatrix}
	      \mblock(0,3)\timestep\beta I_\noutput(3,1)	\mblock(3,3)(4,1)				\mblock(7,3)(3,1)	\mblock(10,3)(3,1)
	      \mblock(0,2)(3,1)						\mblock(3,2)-\alpha_1 D_{\kiter - 1}(4,1)	\mblock(7,2)(3,1)	\mblock(10,2)(3,1)
	      \mblock(0,1)(3,1)						\mblock(3,1)(4,1)				\mblock(7,1)\ddots(1,1)	\mblock(10,1)(3,1)
	      \mblock(0,0)(3,1)						\mblock(3,0)(4,1)				\mblock(7,0)(3,1)	\mblock(10,0)-\alpha_{\ord}D_{\kiter - \ord}(3,1)
	    \end{blockmatrix}$\label{line:BDF7}
	  \State solve ARE \labelcref{eq:ARE} for $L\kiterind$ and $D\kiterind$\label{line:BDF8}
	  \State $\cK\kiterindTbwd = \frac{1}{\sqrt{\weight}}\cB\timek\transpose L\kiterindT D\kiterindT L\kiterind\transpose\cM\kiterind$\label{line:BDF9}
	  \EndFor
\end{algorithmic}
\end{algorithm}

The inputs of \cref{algo:BDF} are the time-dependent matrices \cA\time, \cM\time, $\dot\cM\time$, \cB\time, \cC\time, the time grid \timesteps, 
and the order $\ord \in \left\lbrace1, 2, 3, 4\right\rbrace$ of the BDF method.
Further, the low-rank factors $L_0,\ldots,L_{\ord-1}$, $D_0,\ldots,D_{\ord-1}$ of the initial values $\bX_0,\ldots,\bX_{\ord-1}$ are required as inputs with
sufficient accuracy to obtain the desired order of convergence
\ord.

These initial values could be computed with a different DRE solver of at least
order \ord.
Instead, we compute the initial solution snapshots with sufficiently small time
steps of lower, but successively increasing, order BDF in a wind-up procedure.
This way, the method becomes entirely self-contained.
\begin{algorithm}[t]
	\caption{BDF Start-up Time-steps}%
	\label{algo:startup}
	\hspace*{\algorithmicindent} \textbf{Input} $\timesteps, \ord, \nord$ \\
	\hspace*{\algorithmicindent} \textbf{Output} $\timesteps$
\begin{algorithmic}[1]
	\If{$\ord > 2$}\label{line:startup1}
	\State $\timestep = \ttime_1 - \Starttime,\quad \tilde\timestep = \frac{\timestep}{2^{\nord}}$\label{line:startup2}
	\State $\tilde \ttime_0 = \Starttime,\quad \tilde \ttime_1 = \Starttime + \tilde\timestep$\label{line:startup4}
	\If{$\ord = 3$}\label{line:startup7}
	  \For{$\kiter = 2,\ldots,\nord + 2$}\label{line:startup8}
	    \State $\tilde \ttime\kiterind = \Starttime + 2^{\kiter - 1}\tilde\timestep$\label{line:startup9}
	  \EndFor
	  \State $\timesteps = \{\tilde \ttime_0, \ldots,\tilde \ttime_{\nord + 2}, \ttime_2, \ldots, \ttime_{\ntimestep}\}$\label{line:startup10}
	  \EndIf
	\ElsIf{$\ord = 4$}\label{line:startup11}
	  \For{$\kiter = 1,\ldots,\nord + 2$}\label{line:startup12}
	    \State $\tilde \ttime_{2\kiter} = \Starttime + 2^{\kiter}\tilde\timestep$\label{line:startup13}
	    \State $\tilde \ttime_{2\kiter + 1} = \Starttime + 3\cdot2^{\kiter - 1}\tilde\timestep$\label{line:startup14}
	  \EndFor
	  \State $\timesteps = \{\tilde \ttime_0, \ldots,\tilde \ttime_{2\nord + 5}, \ttime_4, \ldots, \ttime_{\ntimestep}\}$\label{line:startup15}
	\EndIf
\end{algorithmic}
\end{algorithm}
No extra time steps are required for 
$\ord = 1$, where the initial solution $\bX_0$, which is available, is sufficient
and for $\ord = 2$, where the additional initial solution $\bX_1$ can be computed with one BDF 1 step.
But, we generate extra time steps for $\ord \in \{3,4\}$ with a specific wind-up procedure.

We summarize this wind-up procedure in \cref{algo:startup}, which has the inputs \timesteps, \ord, and \nord. 
With these, the smallest additional time step is of length $\tilde\timestep = \frac{\ttime_1 - \Starttime}{2^{\nord}}$.
This is used to compute small BDF steps with step size $\tilde\timestep$ and increasing order from $1$ up to $\ord - 1$.
Then, BDF steps are computed with time step lengths that are doubled \nord{} times until the time steps coincide with the first 
$\Starttime,\ldots,\ttime_{\ord-1}$.
Usually, we choose $\nord = 10$.
Independently of \ntimestep, this results in $10$ extra time steps for BDF 3 and $22$ for BDF 4.
For implementation details of \cref{algo:BDF,algo:startup}, see~\cite[mess\_bdf\_dre.m]{SaaKB-mmess-all-versions}.

In each step of the BDF method, we solve an algebraic Riccati equation (ARE)
\begin{equation}\begin{aligned}\label{eq:ARE}
      0 = 
      \cC\kiterind\transpose\cS\kiterindT\cC\kiterindT
      +(\dot\cM\kiterindT + \cA\kiterindT)\transpose\bX\kiterindT\cM\kiterindT
      +\cM\kiterind\transpose \bX\kiterindT(\dot\cM\kiterindT + \cA\kiterindT)
      -\cM\kiterind\transpose\bX\kiterindT\cB\kiterindT\cB\kiterind\transpose\bX\kiterindT\cM\kiterindT
\end{aligned}\end{equation}
for the low-rank factors 
$L\kiterindT D\kiterindT L\kiterind\transpose \approx \bX\kiterindT$.
The matrices in \cref{eq:ARE} are constructed in \cref{line:BDF3,line:BDF4,line:BDF5,line:BDF6,line:BDF7} of \cref{algo:BDF}
and the coefficients $\alpha_j$ and $\beta$ can be found in,
e.g.,~\cite[Table 5.2]{AscP98}.
Instead of storing the solutions for each time step, \cref{algo:BDF} computes the feedback gain matrices 
$\cK\kiterind \approx \cK\timek, \kiter = 1,\ldots,\ntimestep$ directly, which
further decreases the storage requirements to \(O(mn)\) per time step.
Note that we look at the on-disk storage after completion of the computation.
Intermediately, we may need more memory.
However, we use the Newton alternating-direction implicit (Newton-ADI) method to solve
the ARE \labelcref{eq:ARE}, which can accumulate the feedback gain matrices directly avoiding assembling the low-rank solution factors.
This is a significant performance advantage in terms of the memory requirements
during the computations, especially for large-scale DREs.

There are several new developments in \cref{algo:BDF} compared to~\cite[Algorithm 3.2]{LanMS15}. 
First, the method has the ability to handle non-autonomous DREs with time-dependent matrices and, especially, non-constant $\cM\time$.
Also, the start-up with $\bX_0,\ldots,\bX_\ord$ for order \ord{} with \cref{algo:startup} is independent of any other DRE solver.
Altogether, the extra computational cost of the non-autonomous method are the memory requirements for the time-dependent coefficient matrices
and the matrix-matrix multiplications $\cM\kiterind\transpose L_{j}$ in \cref{line:BDF6}. 
In the autonomous case, only $\hat\cM\transpose L_{\kiter - 1}$ is computed and 
$\hat\cM\transpose L_{\kiter - 2}, \ldots,\hat\cM\transpose L_{\kiter - \ord}$ are reused from the previous time steps
since $\hat\cM$ is constant in time.
Further, \cref{algo:BDF} is embedded in the freely available software package \mmess{} (\cite{SaaKB-mmess-all-versions}), where it
benefits from efficient solvers for the ARE~\labelcref{eq:ARE} and can handle various forms of structures in the coefficient matrices.

\paragraph{Non-autonomous Splitting Schemes.}
The splitting schemes for autonomous DREs utilize the fact that there are closed-form solutions for the two subproblems that arise when treating the affine and the quadratic part of the equation separately. This idea can be partially extended also to the non-autonomous case. 
The relevant subproblems are
\begin{align}
  \label{eq:split_F}
  \cM\transpose \dot{\bX_F} \cM
  &= \cC\transpose\cC
    + (\cA + \dot{\cM})\transpose \bX_F \cM
    + \cM\transpose \bX_F (\cA + \dot{\cM}), \\
  \label{eq:split_G}
  \cM\transpose \dot{\bX_G} \cM
  &= -\cM\transpose\bX_G\cB\cB\transpose\bX_G\cM.
\end{align}
The solution to \cref{eq:split_G} is given implicitly by
\begin{equation}
  \label{eq:G_sol}
  \bX_G\time = {\Big( I + \bX_G(s) \int_{s}^{\ttime}{ \cB(\tau)\cB(\tau)\transpose \diff{\tau}} \Big)}^{-1} \bX_G(s)
\end{equation}
for $s \le \ttime$. This is easily seen by multiplying by the inverted term from the left and then differentiating. 
In the autonomous case (see \cite{Sti18a}), the corresponding formula is 
\begin{equation*}
  \label{eq:G_sol_aut}
  \bX_G\time = {\Big(I + (\ttime -s)\bX_G(s) \cB\cB\transpose \Big)}^{-1} \bX_G(s),
\end{equation*}
and the only difference when solving this subproblem in the non-autonomous case is, thus, that the integral $\int_{s}^{\ttime}{ \cB(\tau)\cB(\tau)\transpose \diff{\tau}}$ needs to be evaluated rather than simply $\cB\cB\transpose$. 

For the affine subproblem, we introduce the two-parameter semigroup
\begin{equation*}
  \cT(r,\ttime) = \texp{\int_{r}^{\ttime}{{\big( (\cA(\tau) + \dot{\cM}(\tau)){\cM(\tau)}^{-1}\big)}\transpose\diff{\tau}}} .
\end{equation*}
Denote the integrand by $\integrand(\tau)$. If $\integrand\time\integrand(s) = \integrand(s)\integrand\time$ for any \ttime{} and $s$, then also the integrals commute:
\begin{equation*}
  \int_{r}^{\ttime}{\integrand(\tau)\diff{\tau}} \int_{r}^{s}{\integrand(\tau)\diff{\tau}} = \int_{r}^{s}{\integrand(\tau)\diff{\tau}} \int_{r}^{\ttime}{\integrand(\tau)\diff{\tau}}.
\end{equation*}
As a consequence, $\big((\cA\time + \dot{\cM}\time)\cM\time^{-1}\big)\transpose$ commutes with $\cT(r,\ttime)$ and
\begin{equation*}
  \frac{\mathrm{d}}{\mathrm{d}\ttime} \cT(r,\ttime)
  ={\big((\cA(\tau) + \dot{\cM}(\tau)){\cM(\tau)}^{-1}\big)}\transpose \cT(r,\ttime).
\end{equation*}
Then, the solution to \cref{eq:split_F} is given implicitly by
\begin{equation}
  \label{eq:sol_F}
    \bX_F\time = \cT(r, \ttime) \bX_F(r) \cT(r,\ttime)\transpose + \int_{r}^{\ttime}{ \cT(s,\ttime) \cM(s)\invtranspose \cC(s)\transpose\cC(s) {\cM(s)}^{-1} \cT(s,\ttime)\transpose \diff{s}},
\end{equation}
for $\Starttime \le r \le \ttime$. This follows quickly by differentiating each of the two terms separately using the above identity and noting that the first term satisfies \cref{eq:split_F} without the $\cC\time\transpose\cC\time$ term, and the second satisfies \cref{eq:split_F} but with $X_F(r) = 0$.

The condition that 
\begin{equation}\label{eq:commcond}
	\integrand\time\integrand(s) = \integrand(s)\integrand\time
\end{equation} 
for any \ttime{} and $s$ does not hold in general, and needs to be verified. Note that it does hold, e.g., when $\cA\time = \alpha\time \bar{\cA}$ and $\cM\time = \mu\time \bar{\cM}$, with
constant matrices $\bar{\cA}$ and $\bar{\cM}$ and scalar functions $\alpha\time$ and $\mu\time$.
This is a common application, e.g., when considering heat flow with a variable thermal conductivity.
In case a different approach can be used to solve \cref{eq:split_F} efficiently, the condition in \cref{eq:commcond} is not required.

In the low-rank setting, we can utilize the Sherman-Morrison-Woodbury matrix inversion lemma 
to express the solution $X_G$
to the nonlinear subproblem \labelcref{eq:split_G} in low-rank form. Given that $\bX_G(s) = L_G D_G L_G\transpose$, from~\Cref{eq:G_sol} we get
\begin{equation}
  \label{eq:split_G_sol_lowrank}
  \bX_G\time = L_G {\Big( I + D_G L_G\transpose \int_{s}^{\ttime}{ \cB(\tau)\cB(\tau)\transpose \diff{\tau}} L_G \Big)}^{-1} D_G L_G\transpose.
\end{equation}
Since the integral does not depend on $\bX_G\time$ and the integrand has a low-rank factorization for each $\tau$, it can easily and efficiently be approximated by applying a quadrature formula followed by a column compression step. This can potentially even be done in an offline phase.

For the affine problem, the constant term in \cref{eq:sol_F} is in low-rank form if the initial condition $\bX_F(r)$ is, and so is the integrand in the integral term. We can, therefore, apply a quadrature rule to approximate the integral and compress the columns of the resulting sum to acquire an approximation to $\bX_F\time$. For each of the terms of this sum, we need to evaluate $\cT(c,\ttime)L$ for different values of $c$ and $L$. For this, we utilize the fact that, at time $s = \ttime$, it is the solution to the system
\begin{equation}
  \label{eq:split_T_system}
  \cM(s)\transpose \dot{Y}(s) = \big( \cA(s) + \dot{\cM}(s) \big)\transpose Y(s), \quad Y(r) = L.
\end{equation}

In the autonomous case, a single integral can be approximated once and then reused in every splitting step, but in the non-autonomous case we need to approximate an integral term in every step. In the autonomous case, we may additionally use the semigroup property $T(r,\ttime)L = T(r,s)T(s, \ttime)L$ to avoid recomputing the same values several times in the quadrature formula. This is no longer possible in the non-autonomous case, because the semigroup will be applied to different $L$ at different quadrature nodes. As a consequence, most of the benefits arising from the splitting, that make these methods very competitive in the autonomous case, are lost in the non-autonomous case.

We summarize the procedure in \cref{algo:splitting}. For clarity, we only consider the Lie splitting case, where we first solve the nonlinear subproblem over a full time step, then the affine subproblem over a full time step. Other splitting schemes will use different combinations of these subproblems, as can be found in~\cite{Sti18a}. However, it is straightforward to adapt \cref{algo:splitting} to those cases. The output of \cref{algo:splitting} is the same as in \cref{algo:BDF}, and so are the basic inputs. The only difference is in the parameters ${\{c_j\}}_{j=1}^s$ and ${\{w_j\}}_{j=1}^s$, which denote the quadrature nodes and weights of a quadrature formula.

\begin{algorithm}[t]
  \caption{Non-autonomous (Lie) splitting scheme}%
  \label{algo:splitting}
	\hspace*{\algorithmicindent} \textbf{Input} $\cA\time,\cM\time,\dot\cM\time,\cB\time,\cC\time,\weight,\timesteps, L_0, D_0$, ${\{c_j\}}_{j=1}^s$, ${\{w_j\}}_{j=1}^s$ \\
	\hspace*{\algorithmicindent} \textbf{Output} $\cK\kiterind, \kiter = 0,\ldots,\ntimestep-1$
\begin{algorithmic}[1]

  \State Invert the direction of time in all matrix functions, e.g., $\cA\time \to \cA(\Endtime  + \Starttime - \ttime)$\label{line:split1}
  \For{$\kiter = 0,\ldots,\ntimestep-1$}\label{line:split2}

    \State Set $\hat{L} = L\kiterind$ and compute $\hat{D} = {\Big( I + D\kiterind L\kiterind\transpose \int_{\ttime\kiterind}^{\ttime\kiternextind}{ \cB(\tau)\cB(\tau)\transpose \diff{\tau}} L\kiterind \Big)}^{-1} D\kiterind$ as in \cref{eq:split_G_sol_lowrank}\label{line:split3} 

    \State Compute $\tilde{L} = \cT(\ttime\kiterind, \ttime\kiternextind) \hat{L}$ by solving \cref{eq:split_T_system} and set $\tilde{D} = \hat{D}$\label{line:split4}

    \State Approximate the integral term in \cref{eq:sol_F} by quadrature $\sum_{j = 1}^{s}{ \tilde{L}_j \tilde{D}_j \tilde{L}_j\transpose }$, where
    $\tilde{L}_j = T(\ttime\kiterind + c_j, \ttime\kiternextind) \cM(\ttime\kiterind+c_j)\invtranspose \cC(\ttime\kiterind+c_j)\transpose$, $\tilde{D}_j = \tau\kiterind w_j I$\label{line:split5}
    
    \State Column compress $L\kiternextind D\kiternextind L\kiternextind\transpose$ with
    $L\kiternextind =
    \begin{bmatrix}
      \tilde{L} & \tilde{L}_1 & \cdots & \tilde{L}_s 
    \end{bmatrix}$
    and
    $D\kiternextind = \text{blkdiag}(\tilde{D}, \tilde{D}_1, \ldots, \tilde{D}_s)$\label{line:split6} 
    
    \State $\cK_{\ntimestep-\kiter-1} = \frac{1}{\sqrt{\weight}}\cB(\ttime\kiternextind)\transpose L\kiternextind D\kiternextind L\kiternextind\transpose\cM(\ttime\kiternextind)$\label{line:split7} 
  \EndFor
\end{algorithmic}
\end{algorithm}

\section{Time stepping Schemes for Closed-loop integrators}\label{sec:time}
Once we have computed a feedback gain \cK\time{} with a method from the previous section,
we can apply it to the solution of the non-linear closed-loop system~\labelcref{eq:closedloop},
as denoted in \cref{eq:feedbackcontrol}, and get a closed-loop input $\controlf\hFEM(\ttime, \sysstateDiff\hFEM\time)$.
In that case, the numerical simulation of this closed-loop 
system can be challenging if $\controlf\hFEM(\ttime, \sysstateDiff\hFEM\time)$ has very large variation 
as we demonstrate in \cref{sec:numeric},
where numerical issues can cause blow-ups.
Similar numerical issues occur in~\cite{FaiW18} for FSI problems with, e.g., the implicit Euler method and the trapezoidal rule.
There, the numerical issues can be overcome with time-adaptive fractional-step-theta schemes.
Since these FSI problems model moving interfaces and use time-dependent meshes,
similar to some of the problems we discuss in \cref{sec:numeric},
we follow~\cite{FaiW18} and the references therein 
and overcome the blow-up behavior in our examples with time-adaptive fractional-step-theta schemes, as well.

To this end, we consider the system to be spatially semi-discretized with finite elements on a mesh of triangular cells 
$\meshcells\time = \{\meshcell\time\}$.
The cells \meshcell\time{} cover the domain \domain{} and are time-dependent in order to treat moving interfaces with time-dependent meshes.
The spatial discretization is a problem-related choice, which we specify 
in \cref{sec:numeric}.
Part of this is the finite element space \FEMspace\hFEM{}, which we use to
define $\sysstateDiff\hFEM \in \FEMspace\hFEM$.
With this, the closed-loop system~\labelcref{eq:closedloop} in semi-discretized form reads
\begin{equation}\label{eq:closedloopsemidiscr}
  \sysstateDiffdot\hFEM\time = \genfuncFEM(\ttime; \sysstateDiff\hFEM\time,\, \controlf\hFEM(\ttime, \sysstateDiff\hFEM\time)),
\end{equation}
where the function $\genfuncFEM\colon\timeinterval\times\mathbb{R}^\nstate\times\mathbb{R}^\nout \to \mathbb{R}^\nstate$ 
is the semi-discretized version of the function \genfunc{}.

For the time discretization, we use \ntimestep{} sub-intervals \timeintervalk{} to partition the time-interval \timeinterval{} 
with the time grid \timestepsfwd:
\begin{equation*}
  \timestepsfwd = \{\ttime_0 ,\ldots, \ttime_{\ntimestep} = \Endtime\},\qquad
  \timestepk = \tk - \tkprev.
\end{equation*}
In contrast to the time discretization $\timesteps$ for the DREs, the time steps
in $\timestepsfwd$ are not necessarily equidistant to allow time-adaptivity,
but $\timesteps \in \timestepsfwd$.
With these time steps, we define discrete approximations $(\sysstatek, \controlk)$ by
\begin{equation*}
  \sysstatek \approx \sysstateDiff\hFEM(\tk),\qquad
  \controlk \approx \controlf\hFEM(\tk, \sysstatek),
\end{equation*}
and the discrete approximation $(\sysstatekref,\, \controlkref)$ of the reference trajectory $(\sysstateref,\,\controlref)$ 
on the corresponding time grid:
\begin{align*}
  \timestepsfwdref &= \{\ttime_0 = \ttime_{0, \text{ref}} ,\ldots, \ttime_{\ntimestep, \text{ref}} = \Endtime\} \in \timestepsfwd.
\end{align*}

For the computation of these discrete approximations, we formulate different time stepping schemes in the same framework
by denoting the parameter $\thetak \in [0, 1]$.
With this, the closed-loop system~\labelcref{eq:closedloopsemidiscr} in discrete form can be approximated by
\begin{equation}\label{eq:closedloopdiscr}
  \frac{\sysstatek - \sysstate\kiterprevind}{\timestepk} = \thetak \genfuncFEM(\sysstatek, \controlk)
				+ (1 - \thetak)\genfuncFEM(\sysstate\kiterprevind, \control\kiterprevind).
\end{equation}       
In case the function $\genfuncFEM$ is nonlinear, \cref{eq:closedloopdiscr} can be solved with, e.g., a Newton method, 
otherwise with a direct linear solver.
As a result, we can obtain various time stepping schemes by choosing different parameters \thetak.
For $\thetak = 0$, \cref{eq:closedloopdiscr} is the explicit Euler scheme and
for $\thetak = 1$, the implicit Euler scheme.
With $\thetak = 0.5$, the scheme corresponds to a trapezoidal rule\footnote{In different PDE contexts, like in~\cite{FaiW18}, the scheme with $\thetak = 0.5$ is also called Crank-Nicolson}.
These are some straight-forward choices for \thetak.
Though, the described framework can be used to formulate more time stepping
schemes as well, where \thetak{} can vary in between time steps,
as for the fractional-step-theta scheme, which is described next.

\subsection{Fractional-step-theta Scheme}\label{sec:FT}
To have a more reliable and robust method regarding the numerical issues that we demonstrate in \cref{sec:numeric},
we combine time-adaptivity with the fractional-step-theta scheme.
It is known to be of second-order accuracy and to be A-stable~\cite{BriGP87}.
Another advantage is that it has little numerical dissipation, which is beneficial for many problems. To formulate this method in our framework, we follow~\cite{Fai17} and define the two parameters
\begin{align*}
  \Theta &= 2 - \sqrt{2} \quad\text{and}\quad \beta = \frac{\Theta - 1}{\Theta - 2} = 1 - \sqrt{\frac{1}{2}}.
\end{align*}
In each time step, we perform three sub-steps, in which we solve \cref{eq:closedloopdiscr} with the parameters
\begin{equation*}
  (\Theta_{\kiter,0},\,\Theta_{\kiter,1},\,\Theta_{\kiter,2}) = (\Theta,\, 1-\Theta,\, \Theta),\qquad
  (\timestep_{\kiter,0},\,\timestep_{\kiter,1},\,\timestep_{\kiter,2}) = 
    (\beta\timestepk,\, (1-2\beta)\timestepk,\, \beta\timestepk).
\end{equation*}
For each sub-step, the initial condition is the solution of the previous sub-step.

Next, we describe the time-adaptive strategies in detail, which we combine with
the fractional-step-theta scheme.

\subsection{Time-adaptive Strategy}
The numerical behavior of the methods for solving the closed-loop system significantly depends on the 
discretization in space and time.
Namely, an equidistant time-discretization with too large time steps (\ntimestep{} too small) can lead to several numerical issues.
For instance, for nonlinear FSI problems, numerical blow-ups 
with the implicit Euler scheme and the trapezoidal rule have been observed in~\cite{FaiW18} 
even though these two methods are A-stable.
Related to this, similar issues can appear when solving the closed-loop system~\labelcref{eq:closedloopdiscr}.
This is potentially be caused, e.g., by small weight parameters \weight{} in the cost functional~\labelcref{eq:cost},
which can cause the feedback control $\controlf\hFEM(\tk, \sysstatek)$ to have very large variation.
Such numerical blow-ups can be prevented by fine time-discretizations.
However, since $\controlf\hFEM(\tk, \sysstatek)$ is unknown a priori, also the required time step size is unknown.
Especially, choosing a fine, equidistant time-discretization with possibly smaller time steps than necessary
can become very computationally costly.

To overcome this difficulty, we refine and coarsen the time steps adaptively in the time-intervals 
that might suffer from numerical blow-ups. 
For this, we use an indicator in order to determine whether the time steps are supposed to be 
larger, smaller, or stay the same.
We denote the indicator as $\indicatork = \indicatorfunc$.
\begin{algorithm}[t]
	\caption{Time-adaptive Step-size Control}%
	\label{algo:timeadap}
	\hspace*{\algorithmicindent} \textbf{Input} $\tk, \tknextref, \sysstatek,
    \controlk, \timestepk, \tol,\gamma,r,\indicatordown,\indicatorup$ \\
    \hspace*{\algorithmicindent} \textbf{Output} $\timestepknext, \retrystep$
\begin{algorithmic}[1]
	\State compute $\indicatork \in \left\lbrace\indicatorke, \indicatorku, \indicatorkdtu\right\rbrace$
	\algorithmiccomment{choose one variant}\label{timeadap1}
	\State $\indicatorkt = {\left(\gammastep\frac{\tol}{\indicatork}\right)}^{\rstep}$\label{timeadap2}
	\State $\retrystep = \indicatorkt < \indicatordown$\label{timeadap3}
	\If{$\indicatordown < \indicatorkt < \indicatorup$\label{timeadap4}}
	  \State $\timestepknext = \timestepk$\label{timeadap5}
	\label{timeadap6}
	\Else
	  \State $\timestepknext = \indicatorkt\timestepk$\label{timeadap7}
	\EndIf\label{timeadap8}
	\If{$\tk + \timestepknext > \tknextref$\label{timeadap9}}\algorithmiccomment{ensure to meet the reference time steps}
	  \State $\timestepknext = \tknextref - \tk$\label{timeadap10}
	\EndIf\label{timeadap11}
\end{algorithmic}
\end{algorithm}

This indicator is used in \cref{algo:timeadap}, the time-adaptive step size control,
which is similar to~\cite[Algorithm 3]{FaiW18}.
From the indicator \indicatork{}, the threshold \indicatorkt{} is computed
with the two parameters $0<\gamma \le 1$ and $r > 0$.
The bounds, $\indicatordown \le 1 \le \indicatorup$, for this threshold 
set the interval in which the current solution
is accepted and the step size stays the same. 
Otherwise, the algorithm decides that, either,
the current solution is considered reasonable ($\retrystep = \text{False}$),
if $\indicatorup < \indicatorkt$, and the step size for the next time step is increased
since we assume that we solved with more accuracy than necessary.
Or, the solution will be recomputed ($\retrystep = \text{True}$),
if $\indicatorkt < \indicatordown$, with the smaller step size \indicatorkt\timestepk.
Further, \cref{algo:timeadap} returns the possibly updated step size \timestepknext.

Additionally, we ensure that the time steps from \timestepsfwdref{} are respected in \cref{timeadap9,timeadap10}.
For these time steps, we have computed \cK\kiterind{} as described in \cref{sec:DRE}.
Otherwise, we use linear interpolation to compute \cK\kiterind{} between time steps from \timestepsfwdref.

A crucial part of this strategy are the 
heuristic indicators \indicatork{} computed in \cref{timeadap1}.
Of these, we define three different variants,
which indicate whether the time step size 
should be changed adaptively. 
One is an error-based indicator, which comes with significantly higher computational costs,
while the other two monitor the feedback control and have no extra computational costs.

\paragraph{Error-based Indicator.}
This indicator is based on a heuristic error estimation. 
First, we compute the solution $(\sysstatekt, \controlkt)$ with one step of the fractional-step-theta scheme with step size \timestepk.
Additionally, we compute $(\sysstatek, \controlk)$ with three steps of step size $\frac{\timestepk}{3}$.
Then, the indicator
\begin{equation*}
  \indicatorke = \nrm{\pickfkt(\sysstatekt) - \pickfkt(\sysstatek)}
\end{equation*}
estimates the error between \sysstatekt{} and \sysstatek.
The function $\pickfkt\colon\mathbb{R}^\nstate \to \mathbb{R}^{\nerror}$, with $\nerror \le \nstate$,
elects certain entries from the solution vector. 
It can be the identity as well $(\nerror = \nstate)$.

This indicator is supposed to refine the time steps especially in time-intervals that might suffer from numerical issues,
which can lead to a blow-up behavior.
Since this behavior strongly depends on the time-discretization, \indicatorke{} is expected to be large in these time-intervals.
However, the computation of \indicatorke{} requires extra computational effort.
Three additional time steps are calculated, which quadruples the cost per time step.
We use the approximate solution with the finer discretization, $(\sysstatek, \controlk)$, to continue the time stepping.

The two control-based indicators, that we describe next, come without the extra computational effort.

\paragraph{Absolute Control-based Indicator.}
We assume that the reference time-discretization \timestepsfwdref{} is chosen such that the reference trajectory snapshots
$(\sysstatekref,\, \controlkref)$ are computed with sufficient accuracy.
Thus, no further refinement is required as long as the feedback control \controlk{} is inactive which is the case
as long as the state \sysstatek{} does not deviate from \sysstatekref.
Therefore, no numerical issues occur during the solution of the closed-loop system as long as the 
feedback control is inactive.
On the other hand, refinement of the time steps might be necessary,
especially, when \controlk{} changes quickly.

The indicator
\begin{equation*}
  \indicatorku = \nrm{\controlk - \controlkprev}
\end{equation*}
monitors the change of the feedback control.
This indicator is computationally cheap to evaluate and tailored to the solution of closed-loop systems, i.e., the feedback control,
which is the focus of this manuscript.

\paragraph{Scaled Control-based Indicator.} Since \indicatorku{} monitors an absolute value, 
this indicator might be sensitive when the time step size \timestepk{} is large and less sensitive
with a small time step size. Consequently, it can be strongly problem dependent.
An alternative is to monitor the scaled change of the feedback control:
\begin{equation*}
  \indicatorkdtu = \nrm{\frac{\controlk - \controlkprev}{\timestepk}}.
\end{equation*}
The indicator \indicatorkdtu{} is an approximation to the norm of the gradient of the feedback control. 
Accordingly, it indicates to refine the time steps if \controlk{} has very large variation as well.
At the same time, it is more robust with respected to the time step size.

We demonstrate the behavior of the time-adaptive strategies with several numerical examples in \cref{sec:fwdexamples}.

\section{Numerical Experiments}\label{sec:numeric}
  In this section, we use three numerical examples
  to put the methods from \cref{sec:DRE,sec:time} to the test
  and demonstrate the interplay between the closed-loop system~\labelcref{eq:closedloop} and the DREs.
  In contrast to the theory in \cref{sec:DRE}, 
  we do not perform the change of variables in time here, such that the DREs are solved backwards in time.
  
  With the first example in \cref{sec:DREnum}, we compare the DRE solvers from \cref{sec:DRE}.
  This example is a small-scale approximation of the Laplacian
  with time-dependent coefficients generated artificially for this purpose.
  Due to the size of $25 \times 25$, it is small enough to compute a reference solution 
  and accurate numerical errors.
  
  The next example models the optimal control
  of heat flow in a 2D cross-section of a steel rail
  with a time-varying thermal conductivity and thus a time-dependent \cA-matrix.
  Also for this example, we can compute a reference solution, since the matrices for this example
  are available in different sizes, starting from dimension $109 \times 109$
  and going up to $\nstate = 79\,841$.
  Consequently, we are able to compare the performance of the DRE solvers in terms of error, runtime, 
  and rank of the computed low-rank factorizations.
  
  The third example is the two-dimensional two-phase Stefan problem, which naturally results in
  large-scale time-dependent matrices due to a moving interface.
  This type of time-dependency in the matrices causes that condition \labelcref{eq:commcond}
  is not fulfilled.
  Consequently, the splitting schemes in their current state are not applicable to this example.
  However, the BDF methods can handle these time-dependent matrices, and we confirm that 
  they still behave similarly to the previous examples.
  While we used the previous examples only to investigate the performance of the DRE solvers,
  for this example we additionally apply several of the computed feedback gain matrices
  in the closed-loop system \labelcref{eq:closedloopsemidiscr}.
  Like this, we compute feedback controls for the Stefan problem.
  With these feedback controls, we solve the non-linear closed-loop system~\labelcref{eq:closedloop} using 
  the methods from \cref{sec:time}.
  In doing so, we put special emphasis on the time-adaptive techniques,
  since those can prevent the solution
  of the closed-loop system from blowing up.
  
  The computations are run on a computer with 2 Intel Xeon Silver 4110 (Skylake) CPUs with 8 cores each
  and 192~GB DDR4 RAM on CentOS Linux release 7.9.2009.
  The used software packages are \fenics{} 2018.1.0, \numpy{} 1.20.3, \scipy{} 1.3.3, and \matlab{} R2017b.
  We provide the codes and data for our experiments
  in~\cite{LQR_Stefan_codes}, including execution logs of the experiments and the results.
  
  \subsection{DRE Solver}\label{sec:DREnum}
	A naive method to solve DREs, is to vectorize them. 
	This means that the DRE is transformed into an equivalent 
	vector-valued non-linear ordinary differential equation (ODE).
	Existing numerical solvers for ODEs can then be applied.
	However, this approach does not exploit the low-rank structure of the solution
	and is unfeasible for large-scale problems. 
	Still, for the sufficiently small examples, we use this method 
	and apply \matlab{}'s \textsf{ode*}-functions
	to generate a high order reference DRE solution
	for error comparisons.
    
  \paragraph{Small-scale Academic Example.}
    With this small example, we intend to analyze the convergence behavior of the DRE solvers.
    Therefore, for this example, we compute a reference solution and compare the numerical errors 
    of the BDF methods and splitting schemes
    as well as the order of convergence.
    In particular, we show how the BDF method (\cref{algo:BDF}) of order 3~and~4
    depends on the parameters for the BDF startup time steps from \cref{algo:startup}.
    Since they would have little meaning in this small-scale example, 
    we do not compare the runtimes of these large-scale methods.

    To generate this example, we use
    a finite-difference approximation of the Laplacian on the unit square ${[0, 1]}^2$
    and assemble the matrices $\hat{\cA} \in \mathbb{R}^{\nstate \times \nstate}, 
    \hat{\cB} \in \mathbb{R}^{\nstate \times 3}, \hat{\cC} \in \mathbb{R}^{1 \times \nstate}$.
    For details about the implementation, see~\cite{LQR_Stefan_codes}.
    Then, these matrices are multiplied with time-dependent scalar functions to get the
    time-dependent matrices
    \begin{align*}
	  \cA\time &= \left(1 + \frac{1}{2}\sin(2\pi\ttime)\right)\hat{\cA}, &
	  \cM\time &= \left(2 + \frac{1}{2}\sin(2\pi\ttime)\right)\cI_\nstate,\\
	  \cB\time &= (3 + \cos(\ttime))\hat{\cB}, &
	  \cC\time &= (1 - \min(\ttime, 1))\hat{\cC}.
    \end{align*}
    Additionally, we choose $\nstate = 25$, $\timeinterval = [0, 0.1]$ and 
    compute a reference approximation $\bX_{\text{ref}}\time$ 
    applying \textsf{ode15s} from \matlab{} after vectorization, as mentioned above.
    With this, in order to assess the quality of our approximate solutions, we
    compare the pointwise relative error:
    \begin{equation*}
	  \errDRE\time = \frac{\nrm[2]{\bX\time - \bX_{\text{ref}}\time}}{\nrm[2]{\bX_{\text{ref}}\time}}.
    \end{equation*}
    
    Following the theory, a BDF method with $\ord \le 6$ stages should converge with order \ord{}, 
    as well~\cite[Section 5.2.3, Example 5.9]{AscP98}.
    As described in \cref{sec:low-rank},
    to achieve these convergence orders, the initialization must be sufficiently accurate.
    For $\ord = 3$~(BDF~3) and $\ord = 4$~(BDF~4), 
    we compute the required initial values with the extra time steps from \cref{algo:startup}. 
    Here, we demonstrate the effect of this on the convergence order.
    However, for $\ord = 1$~(BDF~1) and $\ord = 2$~(BDF~2) no extra time steps are required.
    
    \begin{figure}[t]
	  \centering
      \vspace{0em}%
  \tikzexternalenable%
  \tikzsetnextfilename{smallscale_order_startup}%
  \filemodCmp{tikz/smallscale_order_startup.tikz}{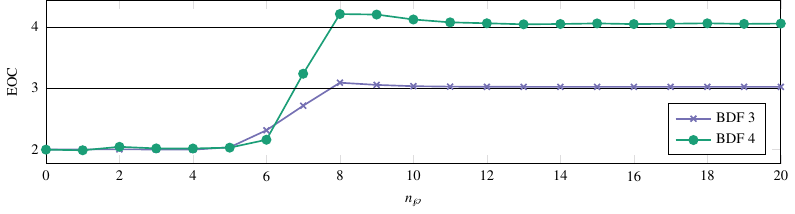}%
  {\tikzset{external/remake next}}{}%
  \begin{tikzpicture}
  \pgfplotstableread{tikz/smallscale/LQR_LTV_startup_x_25_t_512_tend_0.1_rel_err_startup_20_until_9_end_EOC.dat}\tableEOC
  \begin{axis}[
    ylabel={EOC},
    xlabel={\nord},
    xmin=0, xmax=20,
    width=1.04\figurewidth,
    height=.5\figureheight,
    clip=false,
    legend style={at={(0.98,0.05277777)}, anchor=south east},
    ]
      \addlegendentry{BDF 3}
      \addlegendentry{BDF 4}
      
      \addplot [bdf3color, thick, mark=x]
	  table[y index=0, y index=1] {\tableEOC};
	  
      \addplot [bdf4color, thick, mark=*]
	  table[y index=0, y index=2] {\tableEOC};
	  
	  \addplot [black]
	  coordinates {
		  (0,3)
		  (20,3)
		};
		
	  \addplot [black]
	  coordinates {
		  (0,4)
		  (20,4)
		};
  \end{axis}
\end{tikzpicture}%
  \tikzexternaldisable%

      \caption{Convergence orders of BDF with different \nord, small scale example.}%
      \label{fig:orderstartup}
    \end{figure}
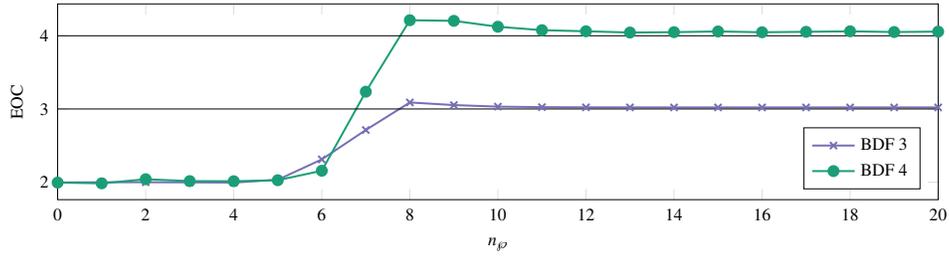
    Correspondingly, \cref{fig:orderstartup} shows the experimental order of convergence (EOC) for BDF~3~and~4.
    In detail, with no extra time steps or \nord{} too small both methods converge with order 2 only.
    At the same time, for $\nord \ge 8$ both methods achieve the expected convergence order \ord{} in this example.
    Moreover, larger values of \nord{} have no (especially also no negative) further effect on the convergence.
    Note, that the choice of \nord{} is problem-dependent, and that for the further computations in this section we use $\nord = 10$.
    
    \begin{figure}[t]
	  \centering
      \vspace{0em}%
  \tikzexternalenable%
  \tikzsetnextfilename{smallscale_order_err}%
  \filemodCmp{tikz/smallscale_order_err.tikz}{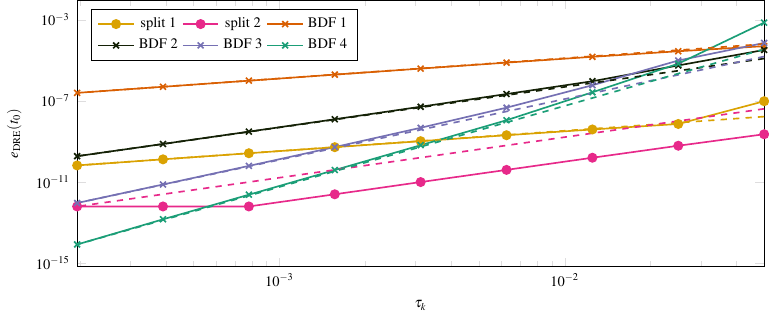}%
  {\tikzset{external/remake next}}{}%
  \begin{tikzpicture}
  \pgfplotstableread{tikz/smallscale/LQR_LTV_compare_x_25_t_512_tend_0.1_rel_err_quad_adaptive_until_9_end.dat}\tableordernine
  \pgfplotstableread{tikz/smallscale/LQR_LTV_compare_x_25_t_512_tend_0.1_rel_err_quad_adaptive_until_9_p_end.dat}\tablepsnine
  
  \begin{axis}[
    ylabel={$\errDRE(\Starttime)$},
    xlabel={\timestepk},
    ymode=log,
    xmode=log,
    xmin=0.00019531, xmax=0.05,
    width=0.98\figurewidth,
    height=.7\figureheight,
    clip=false,
    legend style={at={(0.02,0.96833333)}, anchor=north west},
    legend columns=3,
    name=order,
    ]
      \addlegendentry{split 1}
      \addlegendentry{split 2}
      \addlegendentry{BDF 1}
      \addlegendentry{BDF 2}
      \addlegendentry{BDF 3}
      \addlegendentry{BDF 4}
      \addplot [split1color, thick, mark=*]
	  table[x index=0, y index=1] {\tableordernine};
      \addplot [split2color, thick, mark=*]
	  table[x index=0, y index=2] {\tableordernine};
      \addplot [bdf1color, thick, mark=x]
	  table[x index=0, y index=3] {\tableordernine};
      \addplot [bdf2color, thick, mark=x]
	  table[x index=0, y index=4] {\tableordernine};
      \addplot [bdf3color, thick, mark=x]
	  table[x index=0, y index=5] {\tableordernine};
      \addplot [bdf4color, thick, mark=x]
	  table[x index=0, y index=6] {\tableordernine};
      \addplot [split1color, thick, dashed]
	  table[x index=0, y index=1] {\tablepsnine};
      \addplot [split2color, thick, dashed]
	  table[x index=0, y index=2] {\tablepsnine};
      \addplot [bdf1color, thick, dashed]
	  table[x index=0, y index=3] {\tablepsnine};
      \addplot [bdf2color, thick, dashed]
	  table[x index=0, y index=4] {\tablepsnine};
      \addplot [bdf3color, thick, dashed]
	  table[x index=0, y index=5] {\tablepsnine};
      \addplot [bdf4color, thick, dashed]
	  table[x index=0, y index=6] {\tablepsnine};
  \end{axis}
\end{tikzpicture}%
  \tikzexternaldisable%

      \caption{Convergence of splitting and BDF and the theoretical convergence orders \ord{} (dashed), small scale example.}%
      \label{fig:ordererr}
    \end{figure}
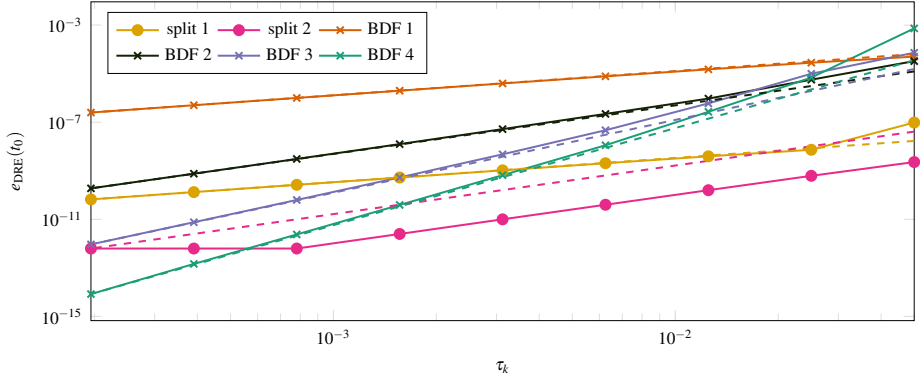
    With this, BDF~1,~2,~3,~and~4 as well as
    the splitting methods of order $\ord = 1$~(split~1) and $\ord = 2$~(split~2)
    can meet their theoretical convergence orders.
    This can be seen in \cref{fig:ordererr},
    which displays the error of the DRE solvers for different time step sizes.
    The observation that the splitting methods yield significantly smaller errors
    than the BDF methods of the same order is also clearly visible here.
    
    \begin{figure}[t]
	  \centering
      \vspace{0em}%
  \tikzexternalenable%
  \tikzsetnextfilename{smallscale_err}%
  \filemodCmp{tikz/smallscale_err.tikz}{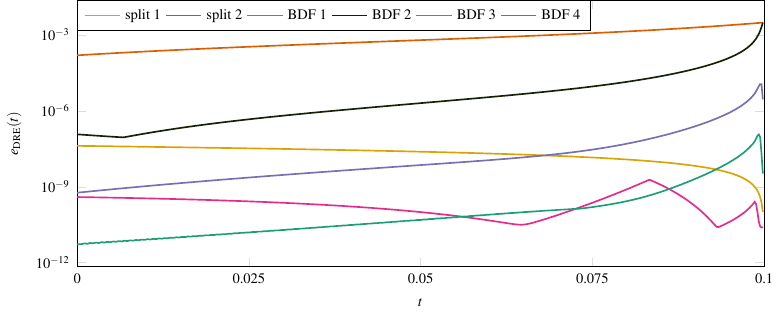}%
  {\tikzset{external/remake next}}{}%
  \begin{tikzpicture}
  \pgfplotstableread{tikz/smallscale/LQR_LTV_compare_x_25_t_512_tend_0.1_rel_err_quad_adaptive_time.dat}\tableerrs
  \begin{axis}[
    ylabel={$\errDRE\time$},
    xlabel={\ttime},
    ymode=log,
    xmin=0, xmax=.1,
    width=0.98\figurewidth,
    height=.7\figureheight,
    xtick = {0,.025,.05,.075,.1},
    xticklabels={0,0.025,0.05,0.075,0.1},
    clip=false,
    legend style={at={(0,1)}, anchor=north west},
    legend columns=6,
    ]
      \addlegendentry{split 1}
      \addlegendentry{split 2}
      \addlegendentry{BDF 1}
      \addlegendentry{BDF 2}
      \addlegendentry{BDF 3}
      \addlegendentry{BDF 4}
      \addplot [split1color, thick]
	  table[x index=0, y index=1] {\tableerrs};
      \addplot [split2color, thick]
	  table[x index=0, y index=2] {\tableerrs};
      \addplot [bdf1color, thick]
	  table[x index=0, y index=3] {\tableerrs};
      \addplot [bdf2color, thick]
	  table[x index=0, y index=4] {\tableerrs};
      \addplot [bdf3color, thick]
	  table[x index=0, y index=5] {\tableerrs};
      \addplot [bdf4color, thick]
	  table[x index=0, y index=6] {\tableerrs};
  \end{axis}
\end{tikzpicture}%
  \tikzexternaldisable%

      \caption{Errors of splitting and BDF for $\ntimestep = 512$ and $\nord = 10$, small scale example.}%
      \label{fig:timeerr}
    \end{figure}
    In addition,
    \cref{fig:timeerr} shows the error of the different DRE solvers
    over the whole time interval for a fixed time step size. 
    The BDF methods start with larger transient errors at \Endtime{}.
    Then, the errors decrease during the backward-in-time solve of the BDF methods.
    Especially BDF~4 is more accurate than the splitting methods in the first half of the time-interval.
    The root cause is that the BDF schemes preserve fixed points, while the splitting schemes do not.
    At the same time, it is important to note that the errors of the splitting methods
    are varying less over the time-interval and stay on roughly the same level.
    Particularly for order $\ord = 1$ and $\ord = 2$,
    the splitting schemes outperform the BDF methods of the same order.
    In essence, the errors in \cref{fig:timeerr} are a typical behavior for these DRE solvers, which 
    is already known for the autonomous variants of BDF and splitting methods~\cite{Sti18a,Lan17}.
    
    For the next example, the BDF and splitting methods benefit from the low-rank factorization of the solution
    since the spatial discretizations can result in large-scale DREs.
    
   \paragraph{Partially Non-autonomous Steel Profile.}
  	The semi-discretized heat transfer model for optimal cooling of steel profiles~\cite{fenicsrail}
  	is an autonomous model.
  	Hence, we modify the example to be non-autonomous and study the performance
  	of the BDF and splitting methods for different matrix sizes and weight parameters
  	in \cref{eq:cost}.
  	We choose this example, because the matrices $\cM$, $\cS$, $\cM_\gamma$, $\cB$ and $\cC$ 
  	are available for different spatial discretizations
  	$\nstate \in \{109, 371, 1\,357, 5\,177, 20\,209, 79\,841\}$ 
  	(see~\cite[\texttt{mess\_get\_linear\_rail.m}]{SaaKB-mmess-all-versions})
  	and we can easily modify the model to be non-autonomous
  	such that it fulfills the condition, which is  necessary for the splitting methods.
  	
  	In detail, the model has $7$ inputs and $6$ outputs
  	and the given matrices form a system equivalent to \cref{eq:sss}:
  	\begin{align*}
		\cM\dot\sysstate\hFEM 	&= (\frac{\kappa}{c\cdot\rho}\cS + \frac{\gamma}{c\cdot\rho}\cM_\gamma)\sysstate\hFEM + \frac{\gamma}{c\cdot\rho}\cB\control\hFEM, \\
		\output\hFEM 			&= \cC\sysstate\hFEM.
	\end{align*}
	Here, the material parameters are the thermal conductivity $\kappa=26.4\frac{kg\,m}{s^3\,^\circ C}$,
	the heat capacity $c=7\,620\frac{m^2}{s^2\,^\circ C}$, 
	the coefficient of thermal conductivity at the input boundary regions 
	$\gamma=7.0164\frac{kg}{s^3\,^\circ C}$, and the density $\rho=654\frac{kg}{m^3}$.
	Following~\cite{BenS05b}, we use the real time-interval $[0, 45s]$,
	which corresponds to model-time $\timeinterval = [0, 4\,500]$.
	Further, to make this system partially non-autonomous, we augment the thermal conductivity
	with a non-linear, smooth, and scalar function:
	\begin{equation*}
		\kappa\time = 26.4 + 0.1 \cdot \left(2 
		+ \cos\left(\frac{2 \cdot \pi \cdot \ttime}{\Endtime}\right)\right)\frac{kg\,m}{s^3\,^\circ C}.
	\end{equation*}
	
  	With this modification the example is similar to, e.g.,~\cite{Lan17}.
  	The difference lies in the scalar function $\kappa\time$, where we
  	specifically picked a non-linear function.
  	Also, due to the constant mass matrix, \cref{algo:BDF} is equivalent to the BDF method therein.
	In addition, the default matrix sizes, number of time steps, 
	and weight parameter for this example
	are given in \cref{tab:pararail}.
	We vary those one by one to examine the behavior of the BDF and splitting methods.
	\begin{table}[t]
	\centering
	\begin{tabular}{|p{4em}|p{4em}|p{4em}|p{4em}|p{4em}|}
		\hline
		\nstate{}	&	\ntimestep{}&	\nin{}	&	\nout{}	&	\weight{}	\\
		\hline
		$109$		&	$128$		&	$7$		&	$6$		&	$1$	 \\
		\hline
	\end{tabular}
	\caption{Partially non-autonomous steel profile parameters (smallest example).}%
	\label{tab:pararail}
	\end{table}
	
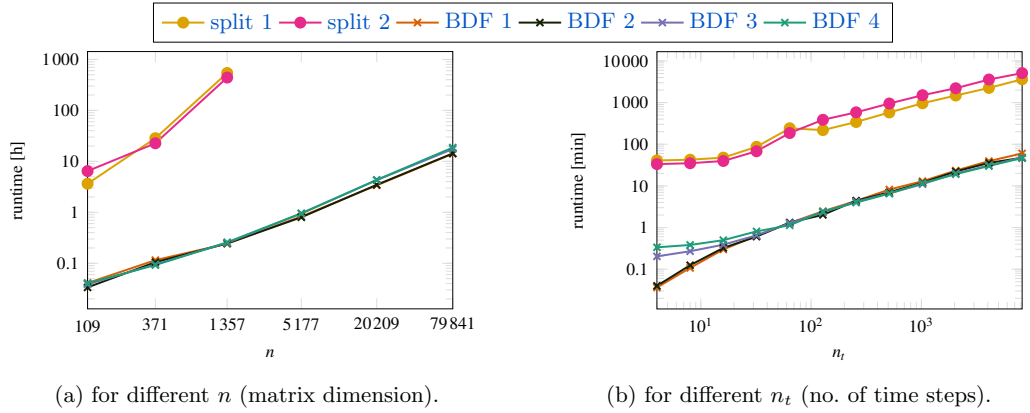
\begin{figure}[t]
    \centering%
    \ref{lname_rail_nx}
    \begin{subfigure}[b]{0.48\textwidth}
    	\centering
  \tikzexternalenable%
  \tikzsetnextfilename{Rail_times_nx}%
  \filemodCmp{tikz/Rail_times_nx.tikz}{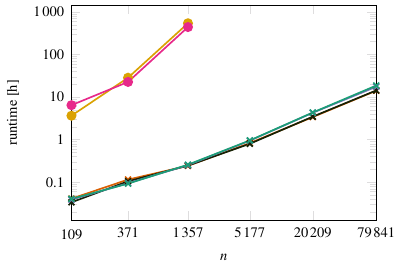}%
  {\tikzset{external/remake next}}{}%
  \begin{tikzpicture}
  \pgfplotstableread{tikz/Rail_A_LTV/A_LTV_plot_data_out.nx_out.time_BDF.dat}\tableBDFtimes
  \pgfplotstableread{tikz/Rail_A_LTV/A_LTV_plot_data_out.nx_out.time_split.dat}\tablesplittimes
  \begin{axis}[
    ylabel={runtime [h]},
    xlabel={\nstate\vphantom{\ntimestep}},
    xmin=109, xmax=79841,
    xmode=log, ymode=log,
    width=0.5\figurewidth,
    height=.6\figureheight,
    clip=true,
    xtick={109,371,1357,5177,20209,79841},
    xticklabels={109\vphantom{$10^1$},371,$1\,357$,$5\,177$,$20\,209$,$79\,841$},
    ytick={0.01,0.1,1,10,100,1000},
    yticklabels={0.01,0.1,1,10,100,$1\,000$},
	legend to name=lname_rail_nx,
    legend columns=-1,
    scaled ticks=false,
    /pgf/number format/1000 sep={\,}
    ]
      \addlegendentry{\small{split 1}}
      \addlegendentry{\small{split 2}}
      \addlegendentry{\small{BDF 1}}
      \addlegendentry{\small{BDF 2}}
      \addlegendentry{\small{BDF 3}}
      \addlegendentry{\small{BDF 4}}

      \addplot [split1color, thick, mark=*]
	  table[x index=0, y index=1] {\tablesplittimes};
	  
      \addplot [split2color, thick, mark=*]
	  table[x index=2, y index=3] {\tablesplittimes};
	  
      \addplot [bdf1color, thick, mark=x]
	  table[x index=0, y index=1] {\tableBDFtimes};
      
      \addplot [bdf2color, thick, mark=x]
	  table[x index=2, y index=3] {\tableBDFtimes};
      
      \addplot [bdf3color, thick, mark=x]
	  table[x index=4, y index=5] {\tableBDFtimes};
      
      \addplot [bdf4color, thick, mark=x]
	  table[x index=6, y index=7] {\tableBDFtimes};
  \end{axis}
\end{tikzpicture}%
  \tikzexternaldisable%

		\caption{for different \nstate{} (matrix dimension).}%
		\label{fig:timeRailnx}
    \end{subfigure}\hfill
    \begin{subfigure}[b]{0.48\textwidth}
    	\centering
  \tikzexternalenable%
  \tikzsetnextfilename{Rail_times_nt}%
  \filemodCmp{tikz/Rail_times_nt.tikz}{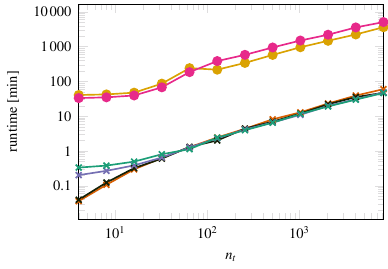}%
  {\tikzset{external/remake next}}{}%
  \begin{tikzpicture}
  \pgfplotstableread{tikz/Rail_A_LTV/A_LTV_plot_data_out.nt_out.time_BDF.dat}\tableBDFtimes
  \pgfplotstableread{tikz/Rail_A_LTV/A_LTV_plot_data_out.nt_out.time_split.dat}\tablesplittimes
  \begin{axis}[
    ylabel={runtime [min]},
    xlabel={\ntimestep},
    xmin=4, xmax=8196,
    xmode=log, ymode=log,
    width=.5\figurewidth,
    height=.6\figureheight,
    clip=true,
    ytick={0.01,0.1,1,10,100,1000,10000},
    yticklabels={0.01,0.1,1,10,100,$1\,000$,$10\,000$},
    /pgf/number format/1000 sep={\,}
    ]

      \addplot [split1color, thick, mark=*]
	  table[x index=0, y index=1] {\tablesplittimes};
	  
      \addplot [split2color, thick, mark=*]
	  table[x index=2, y index=3] {\tablesplittimes};
	  
      \addplot [bdf1color, thick, mark=x]
	  table[x index=0, y index=1] {\tableBDFtimes};
      
      \addplot [bdf2color, thick, mark=x]
	  table[x index=2, y index=3] {\tableBDFtimes};
      
      \addplot [bdf3color, thick, mark=x]
	  table[x index=4, y index=5] {\tableBDFtimes};
      
      \addplot [bdf4color, thick, mark=x]
	  table[x index=6, y index=7] {\tableBDFtimes};
  \end{axis}
\end{tikzpicture}%
  \tikzexternaldisable%

		\caption{for different \ntimestep{} (no.\ of time steps).}%
		\label{fig:timeRailnt}
    \end{subfigure}
    \caption{Runtime of splitting and BDF, steel profile example. 
				All the methods use the same temporal grid.}
\end{figure}
	Starting with different matrix dimensions \nstate{},
	\cref{fig:timeRailnx} displays the corresponding runtimes.
	The first thing to recognize is that the runtime of the BDF methods grows linearly
	with respect to the dimension and 
	there is no significant difference between the orders $\ord \in \{1, \ldots, 4\}$.
	In contrast, the splitting methods require considerably more runtime. 
	As an illustration, while the BDF methods take around $2.3$ minutes
	for the smallest dimension $\nstate = 109$,
	the splitting methods require 3.6~and~6.4~hours. 
	For $\nstate = 371$ it is $6.1$~minutes versus $28.5$~and~$22.6$~hours for BDF and splitting, respectively.
	Forthwith, for dimensions $\nstate > 1\,357$, we stopped the computations
	since the splitting methods did not finish in reasonable time.
	Accordingly, we continue the remaining tests with $\nstate = 109$
	to have reasonable runtimes for the splitting methods 
	even though, like this, the problem is not large-scale.

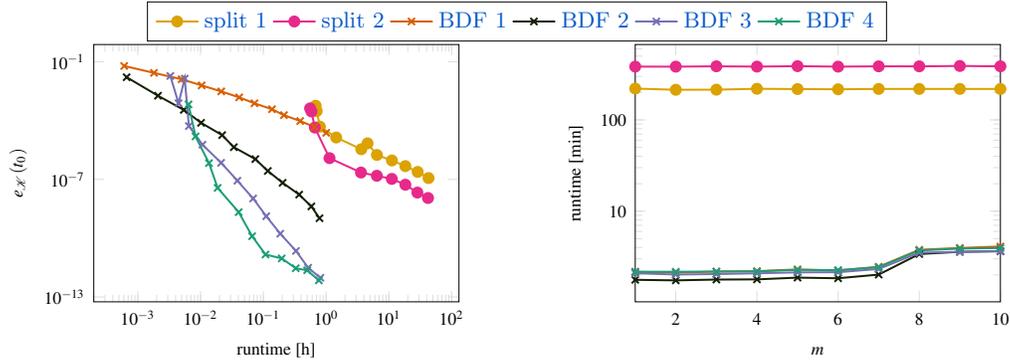
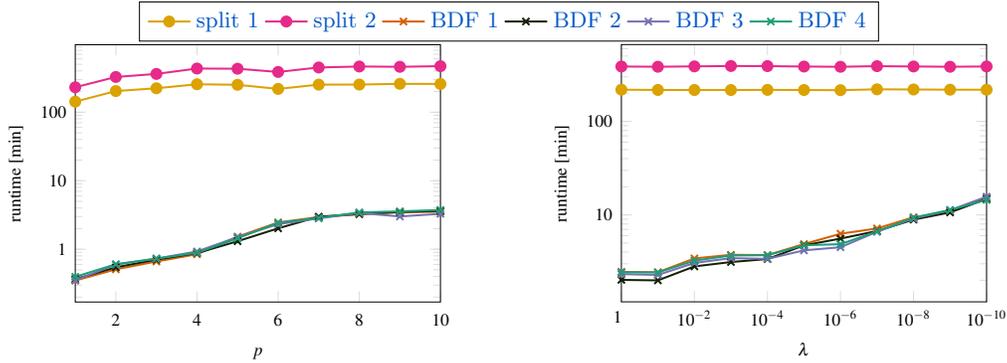
\begin{figure}[t]
    \centering%
    \ref{lname_rail_nx}
    \begin{subfigure}[b]{0.48\textwidth}
    	\centering
  \tikzexternalenable%
  \tikzsetnextfilename{Rail_err_vs_time}%
  \filemodCmp{tikz/Rail_err_vs_time.tikz}{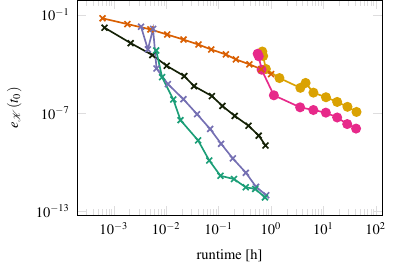}%
  {\tikzset{external/remake next}}{}%
  \begin{tikzpicture}
  \pgfplotstableread{tikz/Rail_A_LTV/A_LTV_plot_err_vs_runtime_ref_diffs_A_LTV_x_109_t_16384_B_7_C_6_1.000000e+00_LTV_fact_1.000000e-01_ode45.mat_BDF.dat}\tableerrvstimebdf
  \pgfplotstableread{tikz/Rail_A_LTV/A_LTV_plot_err_vs_runtime_ref_diffs_A_LTV_x_109_t_16384_B_7_C_6_1.000000e+00_LTV_fact_1.000000e-01_ode45.mat_split.dat}\tableerrvstimesplit
  \begin{axis}[
    ylabel={$e_{\cK}(\Starttime)$},
    xlabel={runtime [h]},
    width=0.5\figurewidth,
    height=.6\figureheight,
    clip=true,
    xmode=log, ymode=log,
    ]

      \addplot [split1color, thick, mark=*]
	  table[x index=0, y index=1] {\tableerrvstimesplit};
	  
      \addplot [split2color, thick, mark=*]
	  table[x index=3, y index=4] {\tableerrvstimesplit};
	  
      \addplot [bdf1color, thick, mark=x]
	  table[x index=0, y index=1] {\tableerrvstimebdf};
      
      \addplot [bdf2color, thick, mark=x]
	  table[x index=3, y index=4] {\tableerrvstimebdf};
      
      \addplot [bdf3color, thick, mark=x]
	  table[x index=6, y index=7] {\tableerrvstimebdf};
      
      \addplot [bdf4color, thick, mark=x]
	  table[x index=9, y index=10] {\tableerrvstimebdf};
  \end{axis}
\end{tikzpicture}%
  \tikzexternaldisable%

		\caption{Runtime vs.\ relative error (compared with \textsf{ode45}),
		$\ntimestep = 16\,384$.}%
		\label{fig:errvstimeRail}
    \end{subfigure}\hfill
    \begin{subfigure}[b]{0.48\textwidth}
    	\centering
  \tikzexternalenable%
  \tikzsetnextfilename{Rail_times_nb}%
  \filemodCmp{tikz/Rail_times_nb.tikz}{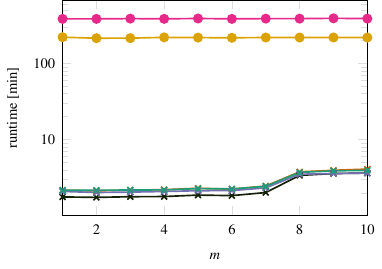}%
  {\tikzset{external/remake next}}{}%
  \begin{tikzpicture}
\tikzstyle{every node}=[font=\small]
  \pgfplotstableread{tikz/Rail_A_LTV/Rail_A_LTV_plot_data_out.nb_out.time_BDF.dat}\tableBDFtimes
  \pgfplotstableread{tikz/Rail_A_LTV/Rail_A_LTV_plot_data_out.nb_out.time_split.dat}\tablesplittimes
  \begin{axis}[
    ylabel={runtime [min]},
    xlabel={\nin \vphantom{[h]}},
    xmin=1, xmax=10,
    ymode=log,
    ytick={0.1,1,10,100,1000},
    yticklabels={0.1,1,10,100,1000},
    xtick={2,4,6,8,10},
    xticklabels={$2\vphantom{10^1}$,$4\vphantom{10^1}$,$6\vphantom{10^1}$,$8\vphantom{10^1}$,$10\vphantom{10^1}$},
    width=.5\figurewidth,
    height=.6\figureheight,
    clip=false,
    ]
	  
      \addplot [split1color, thick, mark=*]
	  table[y index=0, y index=1] {\tablesplittimes};
	  
      \addplot [split2color, thick, mark=*]
	  table[y index=2, y index=3] {\tablesplittimes};
      
      \addplot [bdf1color, thick, mark=x]
	  table[y index=0, y index=1] {\tableBDFtimes};
      
      \addplot [bdf2color, thick, mark=x]
	  table[y index=2, y index=3] {\tableBDFtimes};
      
      \addplot [bdf3color, thick, mark=x]
	  table[y index=4, y index=5] {\tableBDFtimes};
      
      \addplot [bdf4color, thick, mark=x]
	  table[y index=6, y index=7] {\tableBDFtimes};
  \end{axis}
\end{tikzpicture}%
  \tikzexternaldisable%

		\caption{Runtime for different numbers of columns  \nin{} in \cB\time.}%
		\label{fig:timeRailnb}
    \end{subfigure}
    \caption{Runtime of splitting and BDF, steel profile example. 
				All the methods use the same temporal grid.}
\end{figure}

    Next, in \cref{fig:timeRailnt}, for $\ntimestep \in \{2^2, \ldots,
    2^{13}\}$,  we compare the runtime,
    which grows with more time steps, as expected, for all methods.
    Notably, the time per time step is getting slightly smaller
    since the underlying subproblems converge faster.
    A possible explanation is, that in \cref{line:BDF5} of \cref{algo:BDF}
    the quadratic term of the ARE is scaled with $\sqrt{\timestep}$.
    Moreover, the additional overhead to compute the initial values for BDF~3~and~4
    with \cref{algo:startup} is visible for $\ntimestep \in \{2^2, \ldots, 2^{6}\}$.
    
    To better evaluate the efficiency of the methods,
    \cref{fig:errvstimeRail} shows the runtime in relation to the relative error.
    Similar to the small scale-example, we compute a reference solution
    through vectorizing the DRE\@.
    Since not all matrices are time dependent but only $\cA\time$
    in the way that it is a scalar function times a constant matrix
    no stiff ODE solver is required and it is sufficient to use \textsf{ode45}.
Since this example involves large-scale matrices, we do not store the DRE solution
anymore but only the feedback gain matrices $\cK\kiterind$.   
    Thus, for this example, the errors are computed by comparing
    $\cK\kiterind$ at \Starttime{}:
    \begin{equation*}
	  	e_{\cK}(\Starttime) = \frac{\nrm[2]{\cK_0 - \cK_{0,\text{ref}}}}{\nrm[2]{\cK_{0,\text{ref}}}}
    \end{equation*}
    The splitting methods achieve better accuracy than the BDF methods of the same order.
    However, the main computational benefit of the splitting schemes is lost in the non-autonomous case. For an autonomous DRE, only one integral term needs to be approximated for the whole time interval, but now one such integral approximation is needed in each time step. The computational cost is therefore increased by a factor roughly as large as the number of time steps. This is not the case for the BDF schemes, which, therefore, in direct comparison, need significantly less computational effort. The effect can be alleviated by choosing larger tolerances for the splitting subproblem approximations, but as seen in \cref{fig:errvstimeRail}, where this was tried for the 2nd-order splitting scheme, this runs the risk of destroying the overall convergence behavior.
    We also note that the BDF methods converge faster for higher orders, while requiring almost the same runtime for the same \ntimestep{}.
    Thus, BDF~4 is the most efficient method for an accuracy below $10^{-5}$.

\begin{figure}[t]
    \centering%
    \ref{lname_rail_nx}
    \begin{subfigure}[b]{0.48\textwidth}
    	\centering
  \tikzexternalenable%
  \tikzsetnextfilename{Rail_times_nc}%
  \filemodCmp{tikz/Rail_times_nc.tikz}{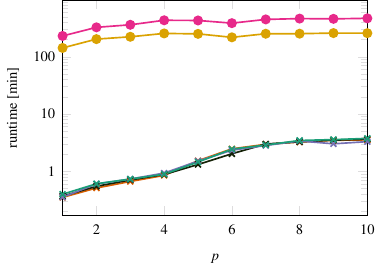}%
  {\tikzset{external/remake next}}{}%
  \begin{tikzpicture}
  \pgfplotstableread{tikz/Rail_A_LTV/Rail_A_LTV_plot_data_out.nc_out.time_BDF.dat}\tableBDFtimes
  \pgfplotstableread{tikz/Rail_A_LTV/Rail_A_LTV_plot_data_out.nc_out.time_split.dat}\tablesplittimes
  \begin{axis}[
    ylabel={runtime [min]},
    xlabel={\nout\vphantom{\weight}},
    xmin=1, xmax=10,
    ymode=log,
    ytick={0.01,0.1,1,10,100,1000,10000},
    yticklabels={0.01,0.1,1,10,100,1000,10000},
    xtick={2,4,6,8,10},
    xticklabels={$2\vphantom{10^1}$,$4\vphantom{10^1}$,$6\vphantom{10^1}$,$8\vphantom{10^1}$,$10\vphantom{10^1}$},
    width=.5\figurewidth,
    height=.6\figureheight,
    clip=false,
    ]
	  
      \addplot [split1color, thick, mark=*]
	  table[y index=0, y index=1] {\tablesplittimes};
	  
      \addplot [split2color, thick, mark=*]
	  table[y index=2, y index=3] {\tablesplittimes};
      
      \addplot [bdf1color, thick, mark=x]
	  table[y index=0, y index=1] {\tableBDFtimes};
      
      \addplot [bdf2color, thick, mark=x]
	  table[y index=2, y index=3] {\tableBDFtimes};
      
      \addplot [bdf3color, thick, mark=x]
	  table[y index=4, y index=5] {\tableBDFtimes};
      
      \addplot [bdf4color, thick, mark=x]
	  table[y index=6, y index=7] {\tableBDFtimes};
  \end{axis}
\end{tikzpicture}%
  \tikzexternaldisable%

		\caption{for different \nout{} (rows in \cC\time).}%
		\label{fig:timeRailnc}
    \end{subfigure}\hfill
    \begin{subfigure}[b]{0.48\textwidth}
    	\centering
  \tikzexternalenable%
  \tikzsetnextfilename{Rail_times_weight}%
  \filemodCmp{tikz/Rail_times_weight.tikz}{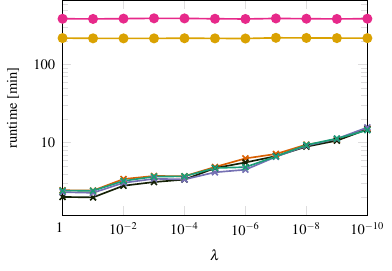}%
  {\tikzset{external/remake next}}{}%
  \begin{tikzpicture}
  \pgfplotstableread{tikz/Rail_A_LTV/A_LTV_plot_data_1_by_out.weight_out.time_BDF.dat}\tableBDFtimes
  \pgfplotstableread{tikz/Rail_A_LTV/A_LTV_plot_data_1_by_out.weight_out.time_split.dat}\tablesplittimes
  \begin{axis}[
    ylabel={runtime [min]},
    xlabel={\weight\vphantom{\nout}},
    xmin=1, xmax=1e+10,
    xmode=log,
    ymode=log,
    xtick={1,100,10000,1000000,100000000,10000000000},
    xticklabels={$1^{\hphantom{1}}$,$10^{-2}$,$10^{-4}$,$10^{-6}$,$10^{-8}$,$10^{-10}$},
    ytick={0.01,0.1,1,10,100,1000,10000},
    yticklabels={0.01,0.1,1,10,100,1000,10000},
    width=.5\figurewidth,
    height=.6\figureheight,
    clip=false,
    ]
	  
      \addplot [split1color, thick, mark=*]
	  table[x index=0, y index=1] {\tablesplittimes};
	  
      \addplot [split2color, thick, mark=*]
	  table[x index=2, y index=3] {\tablesplittimes};
      
      \addplot [bdf1color, thick, mark=x]
	  table[x index=0, y index=1] {\tableBDFtimes};
      
      \addplot [bdf2color, thick, mark=x]
	  table[x index=2, y index=3] {\tableBDFtimes};
      
      \addplot [bdf3color, thick, mark=x]
	  table[x index=4, y index=5] {\tableBDFtimes};
      
      \addplot [bdf4color, thick, mark=x]
	  table[x index=6, y index=7] {\tableBDFtimes};
  \end{axis}
\end{tikzpicture}%
  \tikzexternaldisable%

		\caption{for different \weight{} (weight in cost functional).}%
		\label{fig:timeRailcontrolweights}
    \end{subfigure}
    \caption{Runtime of splitting and BDF, steel profile example. 
				All the methods use the same temporal grid.}
\end{figure}
	We further study the influence of the number of inputs \nin{}, the number of outputs \nout{}, and the weighting factor \weight{} from the performance index~\eqref{eq:cost} on the runtime.
	
    By default, this model has 7 inputs and thus 7 columns in \cB\time{}.
    On the one hand, to test the methods with $\nin \in \{1, \ldots, 7\}$,
    we truncate \cB\time{} after \(m\) columns.
    On the other hand, for $\nin \in \{8, \ldots, 10\}$,
    we generate additional columns that are filled with random values.
    As a result, \cref{fig:timeRailnb} shows that the BDF methods need slightly more time 
    with those random columns, which are not physically motivated in the model.
    Besides this, the runtimes of the BDF and splitting methods are independent of \nin{}.
	
    In contrast, the methods runtimes are affected by different numbers of rows in \cC\time{}.
    \cref{fig:timeRailnc} displays the runtimes for $\nout \in \{1, \ldots, 10\}$.
    Again, we used the first 6 rows in \cC\time{} that are defined in the model. 
    This time, for $\nout \in \{7, \ldots, 10\}$, we add additional rows for physically motivated quantities of interest,
    which are constructed similar to the first 6 outputs.
    The resulting runtimes grow with more rows in \cC\time{},
    but the BDF methods are affected more strongly in this respect.
    This is, because in each BDF step, the underlying iterative ARE solver generates solutions,
    which have \nout{} times the number of iterations of columns.
    Thus, the ARE solution has significantly more columns for larger \nout{} and
    the column compression that we perform on these solutions is significantly more expensive.
    
	 Lastly, the weight \weight{} in the cost functional is an important parameter.
    For instance, the performance of the BDF methods is strongly impacted by \weight{},
    which influences the balance between the quadratic and the constant term
    in the ARE~\labelcref{eq:ARE}.
    In detail, a smaller \weight{} makes the quadratic part more dominant
    and consequently the inner ARE solver in the BDF methods requires more iterations to converge.
    This results in larger runtimes of the BDF methods for smaller \weight{},
    as illustrated in \cref{fig:timeRailcontrolweights}.
    In contrast, the runtimes of the splitting methods are not affected by \weight{} since it only affects \cref{eq:split_G},
    which is solved exactly. However, we note that the temporal discretization errors for both the BDF and splitting schemes depend on the problem being solved and thus on $\weight{}$. To acquire an error of a specific size we might therefore need to use differently sized time steps for different $\weight{}$. The actual computational effort required for a certain accuracy thus varies for both classes of methods. A thorough investigation of these error structures is out of the scope of this paper, but will be the subject of future research.    
    
    Further, the solutions computed by the splitting methods have a smaller numerical rank of $\nLtzero = 68$ 
    compared to $\nLtzero = 84$ for the BDF methods.
    This results in smaller memory requirements for the low-rank solution factors 
    as well as lower computational effort to compute the feedback gain matrices $\cK\kiterind$.
    Ultimately, we apply $\cK\kiterind$ to generate a feedback control
    and do not need the actual DRE solution for this once we have computed $\cK\kiterind$.
    Thus, we can reduce the on-disk memory requirement by storing
    $\cK\kiterind \in \mathbb{R}^{\nstate \times \ninput}$ only,
    which is of fixed size and usually much smaller than the low-rank solution factors.
	  
  \paragraph{Two-dimensional Two-phase Stefan Problem.}
	The Stefan problem, which is treated in~\cite{BarBHetal18a, BarBS22a} as well,
	provides large-scale DREs with time-dependent matrices.
	This time-dependency goes beyond condition~\labelcref{eq:commcond}.
	For this example, we demonstrate the interplay between the feedback gain computation
	and the application of the feedback control during the simulation in~\cref{sec:fwdexamples}.

	The two-dimensional two-phase Stefan problem can model
	solidification and melting of pure materials.
	In our test setup, we consider the time-interval $\timeinterval = [0, 1]$
	and the domain \domain{} is a rectangle
	$[0, 0.5] \times [0, 1] \subset \mathbb{R}^2$
	with the initial interface position at height $0.5$.
	This interface~\Gammaint\time{} splits the domain~\domain{} into two regions 
	that correspond to the solid phase~\domains\time{} and the liquid phase~\domainl\time, 
	as illustrated in \cref{fig:domain}. 
	Because the inner phase-boundary can move, its position is time-dependent.
	Thus, also the two regions \domains\time{} and \domainl\time{} are time-dependent, 
	as well as their corresponding boundary regions.
	They are \Gammain\time, \Gammacool\time{}
	and \GammaN\time{} as denoted in \cref{fig:domain}. 
	The boundary region for the input \Gammain\time{}
	is constant in this representation of the domain. 
	However, we do not restrict ourselves to this case and,
	therefore keep the time-dependence in our notation.
  \tikzexternalenable%
  \tikzsetnextfilename{domain}%
  \filemodCmp{tikz/domain.tikz}{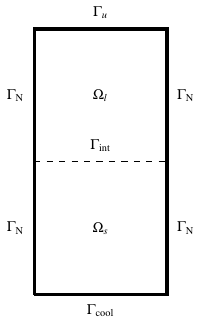}%
  {\tikzset{external/remake next}}{}%
  \begin{figure}[tp]
  \begin{center}
    \begin{tikzpicture}[scale=.45]
      \draw [ultra thick](0,0) -- (5,0) -- (5,10) -- (0, 10) -- (0,0);
      \draw [dashed] (0, 5) -- (5, 5);
      \node at (2.5, -0.6) {$\Gammacool$}; 
      \node at (2.5, 10.6) {$\Gammain$}; 
      \node at (2.5, 7.5) {$\domainl$}; 
      \node at (2.5, 2.5) {$\domains$}; 
      \node at (2.5, 5.6) {$\Gammaint$}; 
      \node at (-0.7,2.5) {$\GammaN$}; 
      \node at (5.7,2.5) {$\GammaN$}; 
      \node at (-0.7,7.5) {$\GammaN$}; 
      \node at (5.7,7.5) {$\GammaN$}; 
    \end{tikzpicture} 
  \end{center}
  \caption{One instance of the domain $\domain \subset \mathbb{R}^2$ of the Stefan problem.} 
  \label{fig:domain}
\end{figure}
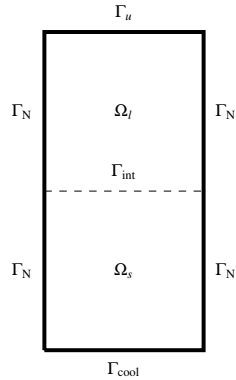%
%
  \tikzexternaldisable%

	We follow the definition of the Stefan problem from~\cite{BarBHetal18a}, but use
	it in a more compact form omitting the couplings with the Navier-Stokes equations
	and the interface graph formulation.
	Both introduce additional nonlinearities to the problem
	and the Navier-Stokes equations add additional DAE structure.
	Before treating these technical challenges, we study the feedback control problem
	in this compact form, without the couplings, to
	develop the general numerical strategy.
	
	The temperature is denoted as $\temp\time$ and modeled by \cref{eq:heat}:
	\begin{subequations}\label{eq:heat}\begin{align}
	  \label{eq:heat1}
	  \timediff\temp - \V \cdot \nabla \temp - \alphac \Delta \temp	&= 0,  																	& &\text{on}\ (0, \Endtime] \times \domain, \\
	  \label{eq:heat2}
	  \partial_{\normal}\temp											&= \control,															& &\text{on}\ (0, \Endtime] \times \Gammain, \\
	  \label{eq:heat3}
	  \temp															&= \Tcool,																& &\text{on}\ (0, \Endtime] \times \Gammacool,\\
	  \label{eq:V1}
	  \Delta \V														&= 0, 																	& &\text{on}\ (0, \Endtime] \times \domain, \\
	  \label{eq:V2}
	  \V 																&= \Big(\frac{1}{\Latent}{[\kk(\nabla \temp)]}_l^s\Big) \cdot \intnormal,	& &\text{on}\ (0, \Endtime] \times \Gammaint.
	\end{align}\end{subequations}
	In \cref{eq:heat1,eq:V2}, the temperature is coupled
	with the extended interface movement~\V\time,
	which we model with the system of algebraic equations~\labelcref{eq:V1}
	and~\labelcref{eq:V2}.
	On the interface, $\restr{\V}{\Gammaint}\time = \Vint\time$
	is the interface movement in the normal direction 
	where \intnormal\time{} is the unit normal vector pointing from \domains\time{}
	to \domainl\time. 
	The extended interface movement~\V\time{} is the smooth extension of
	\Vint\time{} to the whole domain.
	We use \V\time{} to implement
	a spatial discretization that respects the position of \Gammaint\time{}
	and follows the interface movement,
	while preventing mesh tangling and too strong deformations.
	
	\cref{eq:V2} is the Stefan condition, with
	the latent heat constant \Latent{} and the jump of the temperature gradient
	\begin{equation}\label{stefan}
	  {[k(\nabla \temp)]}_l^s = \ks\partial_{\intnormal}\restr{\temp}{\domains} 
	  		- \kl\partial_{-\intnormal}\restr{\temp}{\domainl},
	\end{equation}
	across \Gammaint\time.
	Thus, \V\time{} depends on \temp\time{} as well as
	the heat conductivities in the solid phase \ks{} and in the liquid phase \kl{}
	collected in \alphac\time:
	\begin{equation*}
	  \alphac = \begin{cases}
			    \ks, & \text{on}\ \domains,\\
			    \kl, & \text{on}\ \domainl.\\
			\end{cases}
	\end{equation*}
	Due to this temperature-dependence, the terms $\V \cdot \nabla \temp$ and
	$\alphac \Delta \temp$ in \cref{eq:heat1} are non-linear.
	Further, in \cref{eq:heat2}, we apply the control~\control\time{} as a Neumann condition
	on the control boundary~\Gammain\time.
	Additionally, \cref{eq:heat3} is a Dirichlet condition
	on the cooling boundary~\Gammacool\time{} with the constant \Tcool. 
	A more detailed description including all boundary conditions,
	material parameters, details on the linearization and discretization,
	and the matrix assembly for the DRE can be found in~\cite{BarBS22a}.
	
	In summary, the Stefan problem is a non-linear system of DAEs
	with a moving interface \Gammaint\time{}, 
	which is represented by a moving mesh.
	Our aim is to control the position of \Gammaint\time.
	
	This example results in a large-scale problem with dimension 
	$\nstate = 11\,429$ 
	(max.\ cell width:~0.0268, min.\ cell width:~0.0061, initial mesh).
	For the Stefan problem example, we evaluate the runtime performance of the BDF methods in relation to
	different parameters and matrices \cB\time{} and \cC\time{}. 
	However, the splitting methods are not applicable to this example
	since the time-dependent matrices do not fulfill
	the necessary assumptions in \cref{eq:commcond}.
	To examine their influence on the BDF methods, we individually change 
	the default parameters from \cref{tab:paraStefan}.
	\begin{table}[t]
	\centering
	\begin{tabular}{|p{4em}|p{4em}|p{4em}|p{4em}|p{4em}|}
		\hline
		\ntimestep{}&	\nin{}	&	\nout{}	&	\weight{}	&	\nord{}	\\
		\hline
		$401$		&	$1$		&	$2$		&	$10^{-4}$	&	$10$	 \\
		\hline
	\end{tabular}
	\caption{Two-dimensional two-phase Stefan problem parameters.}%
	\label{tab:paraStefan}
	\end{table}
	
	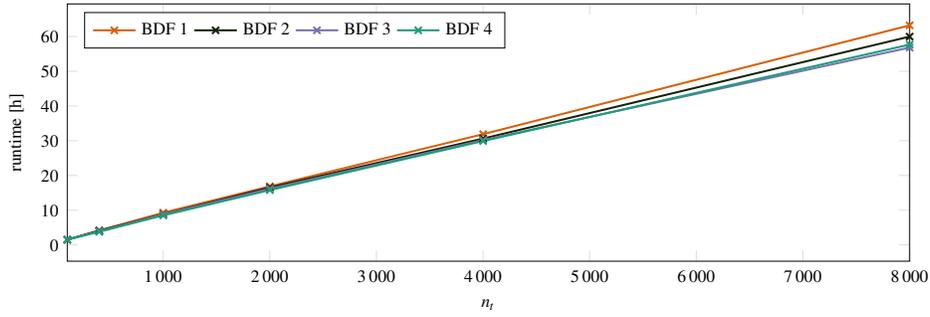
\begin{figure}[t]
		\centering
		\vspace{0em}%
  \tikzexternalenable%
  \tikzsetnextfilename{Stefan_times_nt}%
  \filemodCmp{tikz/Stefan_times_nt.tikz}{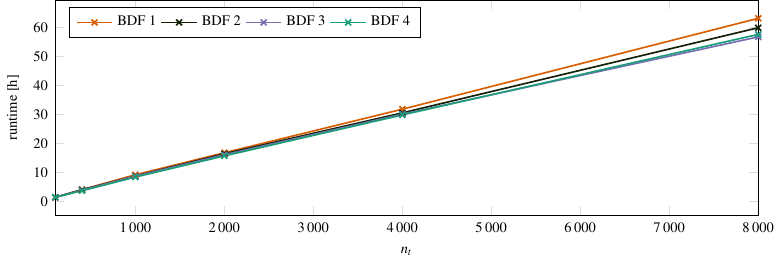}%
  {\tikzset{external/remake next}}{}%
  \begin{tikzpicture}
  \pgfplotstableread{tikz/Stefan_DRE/Stefan_plot_data_out.nt_out.time_1.000000e-04_BDF.dat}\tableBDFtimes
  \begin{axis}[
    ylabel={runtime [h]},
    xlabel={\ntimestep},
    xmin=101, xmax=8001,
    width=\figurewidth,
    height=.6\figureheight,
    clip=true,
    ytick={0,36000,72000,108000,144000,180000,216000},
    yticklabels={0,10,20,30,40,50,60},
    xtick={1000,2000,3000,4000,5000,6000,7000,8000},
    xticklabels={1\,000,2\,000,3\,000,4\,000,5\,000,6\,000,7\,000,8\,000},
    legend style={at={(0.02,0.96833333)}, anchor=north west},
    legend columns=4,
    scaled ticks=false,
    ]
      \addlegendentry{BDF 1}
      \addlegendentry{BDF 2}
      \addlegendentry{BDF 3}
      \addlegendentry{BDF 4}
      
	  \addplot [bdf1color, thick, mark=x]
	  table[x index=0, y index=1] {\tableBDFtimes};
      
      \addplot [bdf2color, thick, mark=x]
	  table[x index=2, y index=3] {\tableBDFtimes};
      
      \addplot [bdf3color, thick, mark=x]
	  table[x index=4, y index=5] {\tableBDFtimes};
      
      \addplot [bdf4color, thick, mark=x]
	  table[x index=6, y index=7] {\tableBDFtimes};
  \end{axis}
\end{tikzpicture}%
  \tikzexternaldisable%

		\caption{Runtime of BDF for different \ntimestep{} (no.\ of time steps),
					Stefan problem example.}%
		\label{fig:timeStefannt}
	\end{figure}
	First, we vary the number of time steps \ntimestep{} from $101$ to $8\,001$. 
	As a result, the runtime of the BDF methods grows linearly with \ntimestep{}
	as can be seen in \cref{fig:timeStefannt}.
	Another key point is that, the runtimes for $\ord = 1, \ldots, 4$ are similar.
	For $\ntimestep = 8\,001$, BDF~3~and~4 are even slightly faster than BDF~1~and~2.
	This results from the fact that the underlying ARE solver (Newton-ADI)
	requires fewer iterations to converge in these cases.
	A possible explanation is that the solutions computed by BDF~3~and~4
	are more accurate resulting in faster convergence.
	
	 Then, as for the example of the steel profile, we vary the number of inputs \nin{},
	the number of outputs \nout{}, and the weight \weight{}
	in the cost functional~\labelcref{eq:cost}.
	The results are consistent with the findings from the previous example 
	and the resulting figures can be found in the supplementary material.
	Moreover, the average numerical ranks of the approximated solutions are
	$\avg(\nLt) = 28$ for BDF~1~and~2, $\avg(\nLt) = 27.9$ for BDF~3, 
	and $\avg(\nLt) = 27.8$ for BDF~4 (for details see supplementary material).
	
	In contrast to the previous example, we evaluate only the low-rank methods 
	since the transformation to an ODE is infeasible to compute.
	Thus, we do not compare any accurate errors in the sense
	that we compare to a reference solution for the DRE\@.
	\begin{figure}[t]
		\centering
		\vspace{0em}%
  \tikzexternalenable%
  \tikzsetnextfilename{Stefan_err_vs_times}%
  \filemodCmp{tikz/Stefan_err_vs_times.tikz}{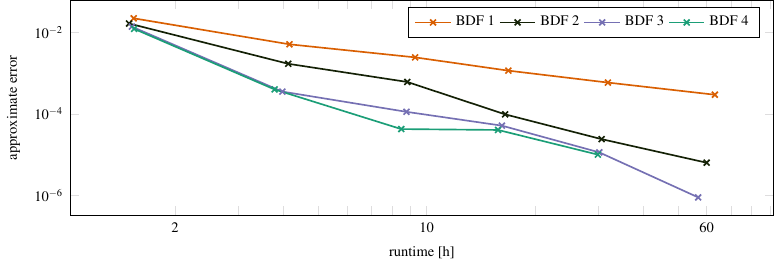}%
  {\tikzset{external/remake next}}{}%
  \begin{tikzpicture}
  \pgfplotstableread{tikz/Stefan_DRE/Stefan_plot_err_vs_time_ref_BDF_4_t_0_nt_8001.dat}\tableerrvstime
  \begin{axis}[
    ylabel={approximate error},
    xlabel={runtime [h]},
    width=\figurewidth,
    height=.6\figureheight,
    clip=true,
    xmode=log, ymode=log,
    xtick={3600,7200,10800,14400,18000,21600,25200,28800,32400,36000,72000,108000,144000,180000,216000,252000,288000,324000,360000},
    xticklabels={1,2,,,,,,,,10,,,,,60,,,,,},
    legend style={at={(0.98,0.96833333)}, anchor=north east},
    legend columns=4,
    ]
      \addlegendentry{BDF 1}
      \addlegendentry{BDF 2}
      \addlegendentry{BDF 3}
      \addlegendentry{BDF 4}

      \addplot [bdf1color, thick, mark=x]
	  table[x index=4, y index=5] {\tableerrvstime};
      
      \addplot [bdf2color, thick, mark=x]
	  table[x index=6, y index=7] {\tableerrvstime};
      
      \addplot [bdf3color, thick, mark=x]
	  table[x index=8, y index=9] {\tableerrvstime};
      
      \addplot [bdf4color, thick, mark=x]
	  table[x index=10, y index=11] {\tableerrvstime};
  \end{axis}
\end{tikzpicture}%
  \tikzexternaldisable%

		\caption{Runtime vs.\ approximate error (compared with BDF~4, $\ntimestep = 8\,001$) of BDF, Stefan problem example. All the methods use the same temporal grid.}%
		\label{fig:errvstimeStefan}
	\end{figure}
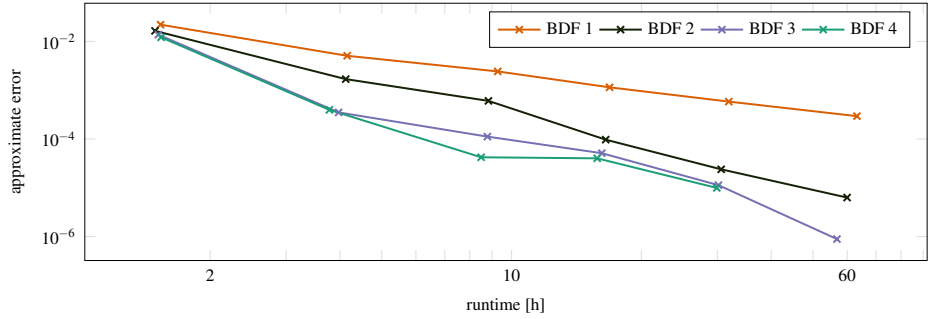
	Instead, we assume the BDF~4 method with the finest time-discretization
	($\ntimestep = 8\,001$) to be most accurate 
	and use that as the reference solution for the comparison
	in \cref{fig:errvstimeStefan}.
	Similar to the results in \cref{fig:errvstimeRail}, 
	the BDF methods are more efficient with higher orders. 
	However, the difference is smaller especially between BDF~3~and~4,
	which might not fully meet their order of convergence.
	
	In the next section, we will use several of the feedback gain matrices that we computed in this paragraph 
	and apply them in the closed-loop system.

  \subsection{Closed-loop System}\label{sec:fwdexamples}
    
	In this section, we use the time stepping methods from \cref{sec:time} 
	to numerically solve the non-linear \cref{eq:closedloop}
	for the two-dimensional two-phase Stefan problem from the previous paragraph.
	In this case, we introduce a perturbation to \cref{eq:heat} 
	by augmenting the Dirichlet boundary condition~\labelcref{eq:heat3} 
	at the bottom of the
	domain with the function \pert\time:
	\begin{equation*}
		    \temp = \Tcool + \pert,	\qquad\text{on}\ (0, \Endtime] \times \Gammacool.
	\end{equation*}
	Here, the control objective is to stabilize the interface position 
	from the discrete reference approximation $(\sysstatekref, \controlkref)$
	and the perturbation \pert\time{} drives the interface away from the
	desired trajectory in the experiments in this section.
	For this, we use a single input ($\nin = 1$) at the top of the domain as depicted in \cref{eq:heat2}.
	Note, that with the present formulation of the Stefan problem in \cref{eq:heat} 
	the interface position can not be measured explicitly.
	Thus, the interface position does not enter the performance index~\labelcref{eq:cost}.
	With this in mind, we use two temperature measurements at the walls \GammaN{} of the domain ($\nout = 2$) as output.
	These two outputs are intended to indicate the deviation of the interface from the reference position.
	
	As a result, we obtain several feedback controls by applying the feedback gain matrices 
	computed in \cref{sec:DREnum} with different weight factors \weight{}.
	
	With these feedback controls, we test the time stepping methods implicit Euler (IE), 
	trapezoidal rule (TR) and our fractional-step-theta scheme (FT) from \cref{sec:FT} 
	to numerically solve \cref{eq:closedloop}.
	For FT, we further compare the different time-adaptivity strategies.
	Those are the error-based indicator (FT-err), the absolute control-based indicator (FT-u) and 
	the scaled control-based indicator (FT-dt-u).
	Additional examples for the Stefan problem 
	with more emphasis on different numbers of inputs 
	and outputs as well as
	different types of perturbations can be found in~\cite{BarBS22a}.
From there, we use the simpler example since the numerical effects that
are highlighted in this manuscript are more clearly visible.
Further, for the more complicated example from~\cite{BarBS22a}, IE as well as TR
are not able to compute a feedback control that successfully stabilizes 
the interface position due to the issues that are discussed next.
	
	    \begin{figure}[t]
		\centering
		\vspace{0em}%
  \tikzexternalenable%
  \tikzsetnextfilename{Stefan_feedback}%
  \filemodCmp{tikz/Stefan_feedback.tikz}{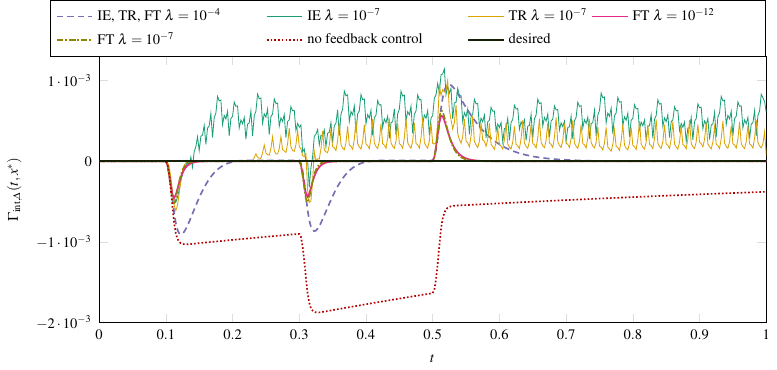}%
  {\tikzset{external/remake next}}{}%
  \begin{tikzpicture}
  \pgfplotstableread{tikz/theta_test/time_401/KK_1.00e-12_1.0e+00_BDF_1_B_1_C_2_1_t_401.mat_tol_u_1.000000e-02_fractional Theta_adaptive_dt_u_diff_not_on_ref.dat}\tableFT
  \pgfplotstableread{tikz/theta_test/time_401/KK_1.00e-04_1.0e+00_BDF_1_B_1_C_2_1_t_401.mat_implicit Euler_diff_not_on_ref.dat}\tableIEtwo
  \pgfplotstableread{tikz/theta_test/time_401/KK_1.00e-07_1.0e+00_BDF_1_B_1_C_2_1_t_401.mat_implicit Euler_diff_not_on_ref.dat}\tableIEthree
  \pgfplotstableread{tikz/theta_test/time_401/KK_1.00e-07_1.0e+00_BDF_1_B_1_C_2_1_t_401.mat_Crank Nicolson_diff_not_on_ref.dat}\tableTRthree
  \pgfplotstableread{tikz/theta_test/time_401/Stefan_test_x_11429_t_401_B_1_C_2_1.000000e-07_BDF_4.mat_tol_u_1.000000e-02_fractional Theta_adaptive_dt_u_diff_not_on_ref_data.dat}\tableFTseven
  \pgfplotstableread{tikz/theta_test/time_401/no_feedback.dat}\tablenofeedback
  \begin{axis}[
    ylabel={$\Gammaintdiff$},
    xlabel={\ttime},
    xmin=0, xmax=1,
    ymin=-2e-3, ymax=1.3e-3,
    width=0.955\figurewidth,
    height=.7\figureheight,
    clip=true,
    restrict y to domain=-.1:.1,
    legend style={at={(1.0,1.0)}, anchor=south east, font=\scriptsize},
    legend cell align={left},
    legend columns=4,
    scaled ticks=false, 
    ]
      \addlegendentry{IE, TR, FT $\weight = 10^{-4}$\hspace{2.3em}\phantom{x}}
      \addlegendentry{IE $\weight = 10^{-7}$\hspace{2.3em}\phantom{x}}
      \addlegendentry{TR $\weight = 10^{-7}$}
      \addlegendentry{FT $\weight = 10^{-12}$\hspace{2.3em}\phantom{x}}
      \addlegendentry{FT $\weight = 10^{-7}$\hspace{2.3em}\phantom{x}}
      \addlegendentry{no feedback control\hspace{2.3em}\phantom{x}}
      \addlegendentry{desired}
      
      \addplot [controled, thick, densely dashed]
	  table[x index=0, y index=1] {\tableIEtwo};
      
      \addplot [extracolor1]
	  table[x index=0, y index=1] {\tableIEthree};
      
      \addplot [extracolor2]
	  table[x index=0, y index=1] {\tableTRthree};
      
      \addplot [extracolor3, thick]
	  table[x index=0, y index=1] {\tableFT};
      
      \addplot [olive, thick, densely dash dot]
	  table[x index=0, y index=1] {\tableFTseven};
      
      \addplot [nofeedbackcolor, thick, densely dotted]
	  table[x index=0, y index=1] {\tablenofeedback};
      \addplot [desired, thick] 
		coordinates {
	  		(0,0)
	  		(1,0)
		};
  \end{axis}
\end{tikzpicture}%
  \tikzexternaldisable%

		\caption{Perturbed interface position relative to the reference trajectory with different feedback controls.}%
		\label{fig:Stefanfeedback}
	    \end{figure}
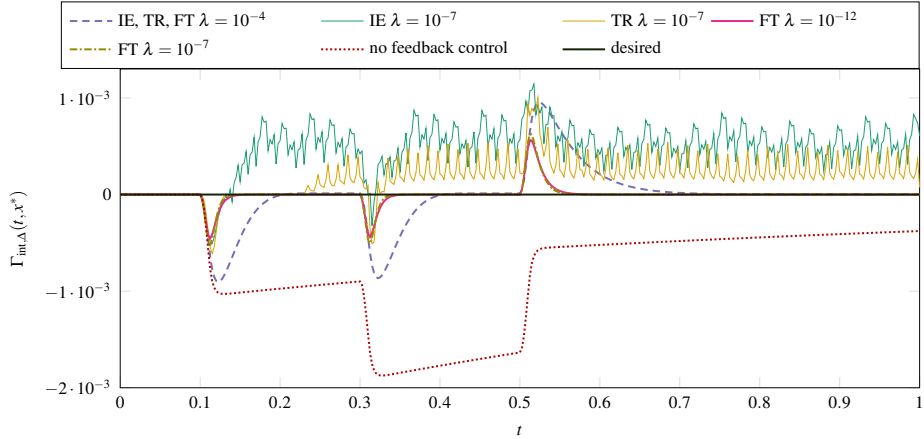
	\cref{fig:Stefanfeedback} displays the relative interface position $\Gammaintdiff = \Gammainttx - \GammaintRef$
	at the point \sysstatexstar{} on the interface which has the largest deviation from the reference trajectory.
	With a weight factor larger than $\weight = 10^{-4}$, 
	the computed feedback control shows very low activity 
	and is not able to steer the interface 
	back towards the reference trajectory in the computed time-frame.
	Decreasing \weight{} means that the performance index~\labelcref{eq:cost} is more dominated by the output deviation
	and the control cost term has less impact. 
	Thus, a smaller \weight{} leads to a more active feedback control, which prevents the interface from
	deviating. 
	In \cref{fig:Stefanfeedback} this is demonstrated with $\weight = 10^{-4}$
	for IE, TR, and FT\@.
	Note that all three lead to very similar results 
	and are summarized in one line for this value of \weight{}.
	
	For smaller \weight{}, the interface can be steered back even faster but IE and TR 
	show a numerical blow-up behavior similar to the experiments in~\cite{FaiW18}.
	This is shown in \cref{fig:Stefanfeedback} for $\weight = 10^{-7}$.
	In our experiments, these blow-ups occur when the feedback control has very large variation and the 
	time step size is too large to properly resolve this.
	
	To overcome this issue, we use time-adaptive FT-dt-u.
	In \cref{fig:Stefanfeedback}, this is demonstrated with $\weight = 10^{-12}$ and FT-dt-u ($\tol = 10^{-2}$).
	Here, the computed feedback control is able to steer the interface position back 
	to the reference trajectory shortly after the perturbation.
	
	The time-adaptive FT-dt-u comes with extra computational cost compared to
	IE and TR, which both require 15.6 minutes to numerically solve
	the closed-loop system with $\ntimestep = 401$ time steps.
	On the other hand, for the same closed-loop system ($\weight = 10^{-4}$), 
	FT-dt-u computes $5\,761$ time steps and requires 546.3 minutes.
	This is partly due to the long time that is needed to steer the interface back to the
	desired trajectory.
	Hence, the control is active in this time period and the time-adaptivity is very expensive.
	In contrast, with $\weight = 10^{-12}$, the interface is back to the reference position in shorter 
	time and the feedback control becomes inactive earlier as well as the time-adaptivity.
	 Thus, FT-dt-u computes $2\,405$ time steps and requires 232.3 minutes.
	However, IE and TR fail with $\ntimestep = 401$ and $\weight \le 10^{-7}$.
The number of time steps $\ntimestep = 401$ is chosen as default
since we observed that this is the smallest \ntimestep{} for which 
the error of the DRE solution and the simulation of the reference trajectory 
stabilize.
	
    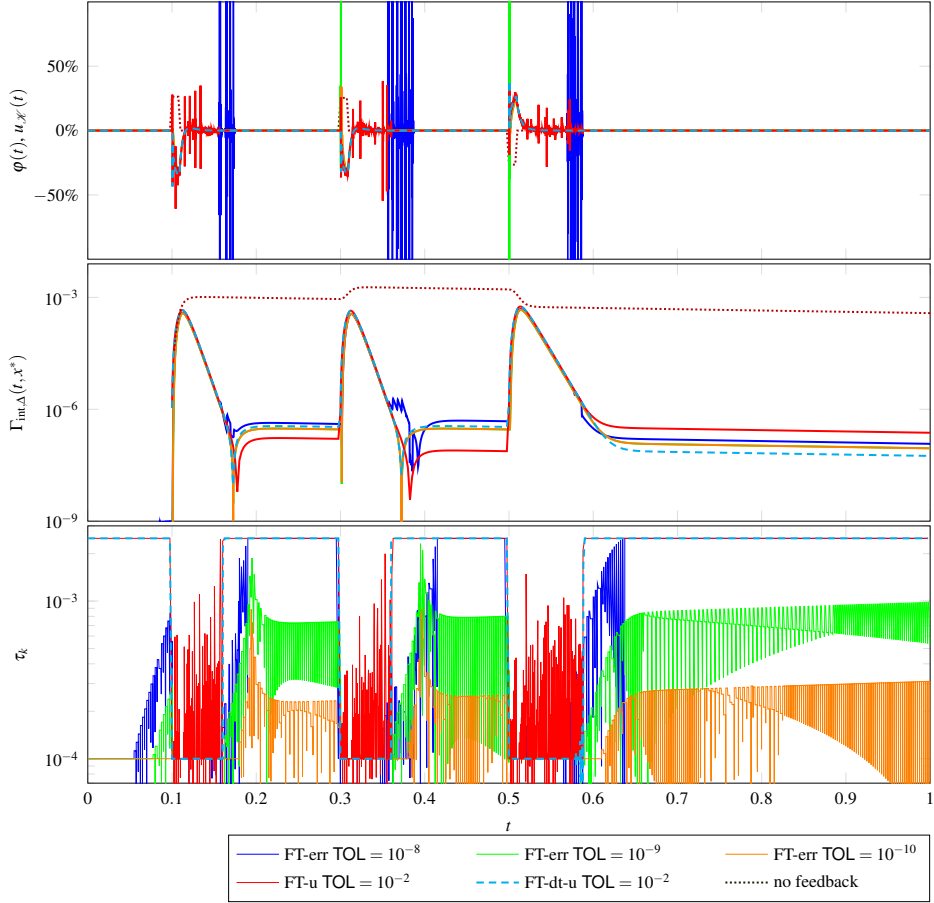
\begin{figure}[tp]
		\centering
        \vspace{0em}%
  \tikzexternalenable%
  \tikzsetnextfilename{Stefan_feedback_adaptive}%
  \filemodCmp{tikz/Stefan_feedback_adaptive.tikz}{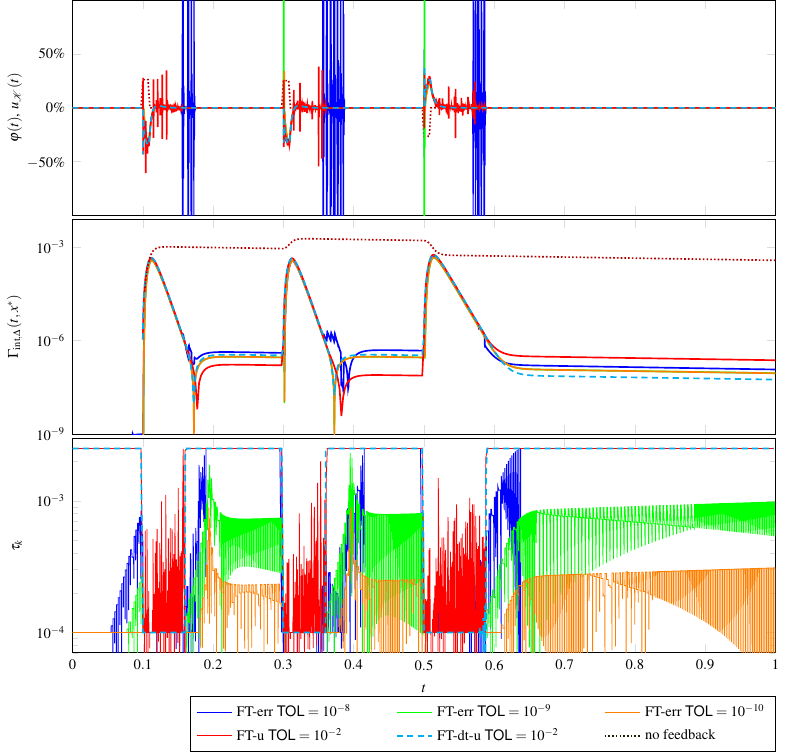}%
  {\tikzset{external/remake next}}{}%
  \begin{tikzpicture}
  
  \pgfplotstableread{tikz/theta_test/time_401/KK_1.00e-12_1.0e+00_BDF_1_B_1_C_2_1_t_401.mat_tol_u_1.000000e-02_fractional Theta_adaptive_dt_u_diff_not_on_ref_full.dat}\tabledtu
  \pgfplotstableread{tikz/theta_test/time_401/KK_1.00e-12_1.0e+00_BDF_1_B_1_C_2_1_t_401.mat_tol_1.000000e-02_fractional Theta_adaptive_u_diff_not_on_ref_full.dat}\tableu
  \pgfplotstableread{tikz/theta_test/time_401/KK_1.00e-12_1.0e+00_BDF_1_B_1_C_2_1_t_401.mat_tol_err_1.000000e-08_fractional Theta_adaptive_error_diff_not_on_ref_full.dat}\tableerreight
  \pgfplotstableread{tikz/theta_test/time_401/KK_1.00e-12_1.0e+00_BDF_1_B_1_C_2_1_t_401.mat_tol_err_1.000000e-09_fractional Theta_adaptive_error_diff_not_on_ref_full.dat}\tableerrnine
  \pgfplotstableread{tikz/theta_test/time_401/KK_1.00e-12_1.0e+00_BDF_1_B_1_C_2_1_t_401.mat_tol_err_1.000000e-10_fractional Theta_adaptive_error_diff_not_on_ref_full.dat}\tableerrten
  \pgfplotstableread{tikz/theta_test/time_401/no_feedback.dat}\tablenofeedback
  \begin{axis}[
    ylabel={\pert\time, \controlf\time},
	ytick={},
    xmin=0, xmax=1,
    ymin=-1, ymax=1,
    xtick = {0,0.1,0.2,0.3,0.4,0.5,0.6,0.7,0.8,0.9,1},
    xticklabels ={,,,,,,,,,,},
    ytick = {-0.5,0,0.5},
	yticklabels={\phantom{$10^{-3}$}\makebox[0pt][r]{$-50\%$},$0\%$,$50\%$},
    width=\figurewidth,
    height=.6\figureheight,
    clip=true,
    name=feedbackadaptivecontrol,
    ]
      
      \addplot [nofeedbackcolor, thick, densely dotted]
	  table[x index=0, y index=2] {\tablenofeedback};
      
      \addplot [fterr8color]
	  table[x index=0, y index=2] {\tableerreight};
      
      \addplot [fterr9color, thick]
	  table[x index=0, y index=2] {\tableerrnine};
      
      \addplot [fterr10color, thick]
	  table[x index=0, y index=2] {\tableerrten};
      
      \addplot [ftu2color, thick]
	  table[x index=0, y index=2] {\tableu};
      
      \addplot [ftdtu2color, thick, densely dashed]
	  table[x index=0, y index=2] {\tabledtu};
  \end{axis}

  \pgfplotstableread{tikz/theta_test/time_401/KK_1.00e-12_1.0e+00_BDF_1_B_1_C_2_1_t_401.mat_tol_u_1.000000e-02_fractional Theta_adaptive_dt_u_diff_not_on_ref_abs_full.dat}\tabledtuabs
  \pgfplotstableread{tikz/theta_test/time_401/KK_1.00e-12_1.0e+00_BDF_1_B_1_C_2_1_t_401.mat_tol_1.000000e-02_fractional Theta_adaptive_u_diff_not_on_ref_abs_full.dat}\tableuabs
  \pgfplotstableread{tikz/theta_test/time_401/KK_1.00e-12_1.0e+00_BDF_1_B_1_C_2_1_t_401.mat_tol_err_1.000000e-08_fractional Theta_adaptive_error_diff_not_on_ref_abs_full.dat}\tableerreightabs
  \pgfplotstableread{tikz/theta_test/time_401/KK_1.00e-12_1.0e+00_BDF_1_B_1_C_2_1_t_401.mat_tol_err_1.000000e-09_fractional Theta_adaptive_error_diff_not_on_ref_abs_full.dat}\tableerrnineabs
  \pgfplotstableread{tikz/theta_test/time_401/KK_1.00e-12_1.0e+00_BDF_1_B_1_C_2_1_t_401.mat_tol_err_1.000000e-10_fractional Theta_adaptive_error_diff_not_on_ref_abs_full.dat}\tableerrtenabs
  \pgfplotstableread{tikz/theta_test/time_401/no_feedback_abs.dat}\tablenofeedbackabs
  \begin{axis}[
    at={(feedbackadaptivecontrol.below south west)},
    anchor = north west,
    yshift=5pt,
    ylabel={$\Gammaintdiff$},
	ytick={},
    xmin=0, xmax=1,
    ymin=1e-9,
    xtick = {0,0.1,0.2,0.3,0.4,0.5,0.6,0.7,0.8,0.9,1},
    xticklabels ={,,,,,,,,,,},
    ymode=log,
    width=\figurewidth,
    height=.6\figureheight,
    clip=true,
    legend style={at={(1,0)}, anchor=north east},
    name=feedbackadaptive,
    ]
      
      \addplot [fterr8color, thick]
	  table[x index=0, y index=1] {\tableerreightabs};
      
      \addplot [fterr9color, thick]
	  table[x index=0, y index=1] {\tableerrnineabs};
      
      \addplot [fterr10color, thick]
	  table[x index=0, y index=1] {\tableerrtenabs};
      
      \addplot [ftu2color, thick]
	  table[x index=0, y index=1] {\tableuabs};
      
      \addplot [ftdtu2color, thick, densely dashed]
	  table[x index=0, y index=1] {\tabledtuabs};
      
      \addplot [nofeedbackcolor, thick, densely dotted]
	  table[x index=0, y index=1] {\tablenofeedbackabs};
  \end{axis}

  \pgfplotstableread{tikz/theta_test/time_401/KK_1.00e-12_1.0e+00_BDF_1_B_1_C_2_1_t_401.mat_tol_u_1.000000e-02_fractional Theta_adaptive_dt_u_diff_not_on_ref_time_steps_full.dat}\tabledtusteps
  \pgfplotstableread{tikz/theta_test/time_401/KK_1.00e-12_1.0e+00_BDF_1_B_1_C_2_1_t_401.mat_tol_1.000000e-02_fractional Theta_adaptive_u_diff_not_on_ref_time_steps_full.dat}\tableusteps
  \pgfplotstableread{tikz/theta_test/time_401/KK_1.00e-12_1.0e+00_BDF_1_B_1_C_2_1_t_401.mat_tol_err_1.000000e-08_fractional Theta_adaptive_error_diff_not_on_ref_time_steps_full.dat}\tableerreightsteps
  \pgfplotstableread{tikz/theta_test/time_401/KK_1.00e-12_1.0e+00_BDF_1_B_1_C_2_1_t_401.mat_tol_err_1.000000e-09_fractional Theta_adaptive_error_diff_not_on_ref_time_steps_full.dat}\tableerrninesteps
  \pgfplotstableread{tikz/theta_test/time_401/KK_1.00e-12_1.0e+00_BDF_1_B_1_C_2_1_t_401.mat_tol_err_1.000000e-10_fractional Theta_adaptive_error_diff_not_on_ref_time_steps_full.dat}\tableerrtensteps
  \begin{axis}[
    at={(feedbackadaptive.below south west)},
    anchor = north west,
    yshift=5pt,
    ylabel={\timestepk},
    xlabel={\ttime},
    xmin=0, xmax=1,
    ymin=0.00007, ymax=0.003,
    ymode=log,
    width=\figurewidth,
    height=.6\figureheight,
    clip=true,
    legend style={at={(1,-0.2)}, anchor=north east,font= \scriptsize},
    legend cell align={left},
    legend columns=3,
    ]
      \addlegendentry{FT-err $\tol = 10^{-8}$\hspace{2.3em}\phantom{x}}
      \addlegendentry{FT-err $\tol = 10^{-9}$\hspace{2.3em}\phantom{x}}
      \addlegendentry{FT-err $\tol = 10^{-10}$}
      \addlegendentry{FT-u $\tol = 10^{-2}$\hspace{2.3em}\phantom{x}}
      \addlegendentry{FT-dt-u $\tol = 10^{-2}$\hspace{2.3em}\phantom{x}}
      \addlegendentry{no feedback}
      
      \addplot [fterr8color]
	  table[x index=0, y index=1] {\tableerreightsteps};
      
      \addplot [fterr9color]
	  table[x index=0, y index=1] {\tableerrninesteps};
      
      \addplot [fterr10color]
	  table[x index=0, y index=1] {\tableerrtensteps};
      
      \addplot [ftu2color]
	  table[x index=0, y index=1] {\tableusteps};
      
      \addplot [ftdtu2color, thick, densely dashed]
	  table[x index=0, y index=1] {\tabledtusteps};
	  
	  \addplot [desired, thick, densely dotted] 
	coordinates {
	  (0,0.001)
	};
  \end{axis}
\end{tikzpicture}%
  \tikzexternaldisable%

		\caption{Perturbation and different time-adaptive feedback controls 
				(clipped, max.\ value: 8.4, min.\ value: -39.3) (top),
				Perturbed interface position (center), and time step sizes 
				(bottom) for $\weight = 10^{-12}$.}%
	\label{fig:Stefanfeedbackadaptive}
    \end{figure} 
	The corresponding feedback control that is computed with FT-dt-u for $\weight = 10^{-12}$ is shown
	in comparison with  
	the feedback controls computed with the different time-adaptive strategies 
	FT-u and FT-err for the same closed-loop system in \cref{fig:Stefanfeedbackadaptive} (top).
	Further, the figure displays the perturbation \pert\time{} as a dotted line.
	To clarify the y-axis labeling, it is important to note 
	that it is relative to the original boundary value \Tcool{}.
	
	An important parameter for the adaptive step size computation in \cref{algo:timeadap} is \tol.
	With smaller \tol{}, \cref{algo:timeadap} computes smaller time step sizes when the indicator grows.
	On the one hand, a larger \tol{} can save some time steps by increasing the time step size earlier with dropping indicator values and reduce the computational effort.
	On the other hand, \cref{algo:timeadap} can fail to prevent a numerical blow-up
	if \tol{} is chosen too large and the time step sizes grow too early.
	This effect is presented in \cref{fig:Stefanfeedbackadaptive} for FT-u with $\tol = 10^{-2}$ and FT-err with $\tol = 10^{-8}$.
	Both increase the time step size when the feedback control is nearly inactive but still active enough
	to lead to a blow-up. 
	This blow-up then leads to an increase of the indicator, and consequently a decrease of the time step size 
	such that the feedback controls become inactive again.
	The result is that the feedback control is quickly bouncing back and forth between large and small values.
	
	In contrast to this, for FT-dt-u, the indicator monitors the relative change of the feedback control.
	Thus, as long as the feedback control is active the algorithm sets \timestepk{} 
	to the minimum value ($\timestepk = 0.0001$).
	However, when the control is inactive the algorithm sets \timestepk{} back
	to the maximum value ($\timestepk = 0.0025$) as can be seen in \cref{fig:Stefanfeedbackadaptive} 
	(bottom, dashed line).
	We expect this to be the general behavior of FT-dt-u if \tol{} is chosen sufficiently small.
	
	To emphasize a meaningful comparison, for FT-u, the indicator monitors the absolute change of the feedback control.
	This means, with a smaller \tol{}, it behaves the same as FT-dt-u ($\tol < 10^{-5}$).
	In particular, it depends on the magnitude of the feedback control values, and 
	the right choice of \tol{} is therefore strongly problem dependent.
	From our numerical experiments, we expect FT-dt-u to be more robust with respect to the magnitude of the feedback control values 
	and, thus, the problem dependent choice of \tol.
	Further, we can choose a smaller \tol{} to be on the safe side and get the same \ntimestep{} 
	since it jumps back and forth from the maximum to the minimum \timestepk.
	
	Different from the other two, the indicator FT-err monitors the difference of the computed solutions for two different time step sizes.
	This results in a larger computational effort for a single time step compared to the alternative indicators.
	For instance, if the feedback control has very large variation and consequently the solution changes more, then \timestepk{} is reduced.
	However, FT-err can not distinguish between high activity of the reference solution and of the perturbed solution. 
	In addition, it enlarges \timestepk{} gradually.
	Altogether, it requires significantly smaller \tol, larger \ntimestep, and, consequently, increases the runtime to prevent a blow-up.
	To illustrate this, the computed time steps for three different choices of \tol{} are presented in \cref{fig:Stefanfeedbackadaptive} (bottom).
	To clarify the behavior of FT-err, \cref{algo:timeadap} is restricted to meet the reference time steps \timestepsfwdref{} for all the time-adaptive strategies.
	This can result in time step sizes that are smaller than the minimum \timestepk{} 
	if the adaptively computed time steps would otherwise jump past the next time step from \timestepsfwdref.
	Some quickly alternating time step size adjustments computed by FT-err result from this restriction.
    
\begin{figure}[tp]
		\centering
		\vspace{0em}%
  \tikzexternalenable%
  \tikzsetnextfilename{FT_adaptive_times}%
  \filemodCmp{tikz/FT_adaptive_times.tikz}{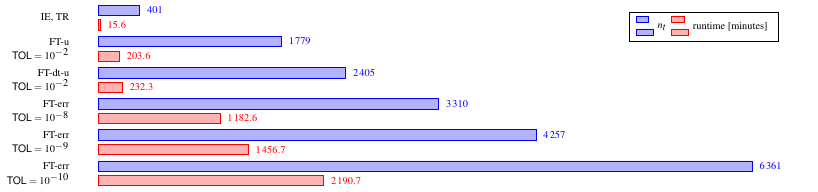}%
  {\tikzset{external/remake next}}{}%
  \begin{tikzpicture}
\tikzstyle{every node}=[font=\tiny]
	\begin{axis}[
	    width=1.1\figurewidth,
	    height=1.2\figureheight,
	    xbar=-2pt,
	    y=15pt,
	    bar width=-5pt,
	    y axis line style = { opacity = 0 },
	    axis x line       = none,
	    legend style={at={(0.5,-0.15)},
	      anchor=north,legend columns=-1},
	    xlabel=x-label,
	    ytick=data,
	    symbolic y coords={IE/TR,FT-u,FT-dt-u,FT-err8,FT-err9,FT-err10},
	    y dir=reverse,
	    yticklabels={{IE, TR},{FT-dt-u\\ $\tol=10^{-2}$},{FT-u\\ $\tol=10^{-2}$},{FT-err\\ $\tol=10^{-8}$},{FT-err\\ $\tol=10^{-9}$},{FT-err\\ $\tol=10^{-10}$}},
	    yticklabel style={align=right,xshift=2em},
	    tickwidth         = 0pt,
	    nodes near coords,
	    legend style={at={(0.95,0.95)}, anchor=north east},
	    /pgf/number format/1000 sep={\,}
	    ]
		\addplot coordinates {(401,IE/TR)	(2405,FT-dt-u)	(1779,FT-u)		(3310,FT-err8)		(4257,FT-err9)		(6361,FT-err10)};
		\addplot coordinates {(15.6,IE/TR)	(232.3,FT-dt-u)	(203.6,FT-u)	(1182.6,FT-err8)	(1456.7,FT-err9)	(2190.7,FT-err10)};
		\legend{\ntimestep,runtime [minutes]}
	\end{axis}
\end{tikzpicture}%
  \tikzexternaldisable%

		\caption{Number of time steps \ntimestep{} and runtimes of IE, TR, 
			and FT with different time-adaptive strategies, $\weight = 10^{-12}$.}%
		\label{fig:Stefanfeedbackadaptivetimesandsteps}
  \end{figure}
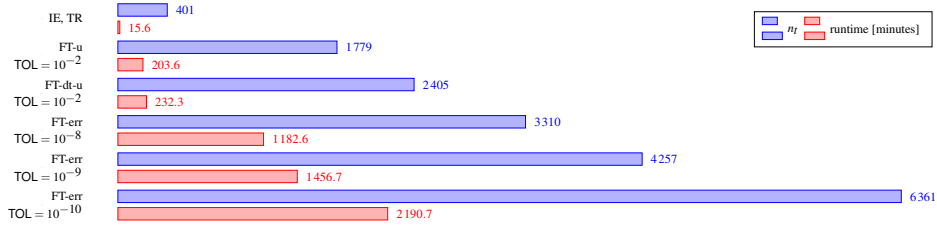

	To illustrate the just described effects and differences between the time-adaptive strategies, 
	\cref{fig:Stefanfeedbackadaptivetimesandsteps} displays the different \ntimestep{} and runtimes.
	Notably, IE and TR require $2.3$ seconds per time step.
	With the latter two methods, it would take 390 minutes to solve the closed-loop system with $10\,000$ equidistant time steps.
	This number of time steps corresponds to the step size that FT-dt-u uses while the control is active.
	Compared to that, FT-u and FT-dt-u require shorter runtime and significantly fewer time steps while using the same small time step size locally.
	Further, they simulate the closed-loop system more reliably in presence of strongly varying controls due to their adaptivity, which is tuned to this problem.
	In detail, the runtime per time step is $6.9$ seconds for FT-u and $5.8$ seconds for FT-dt-u on average.
	FT-u takes longer on average since it discards time steps more often than FT-dt-u since it tries to take larger time steps too early.
	A discarded time step results in a re-computation of this time step with a smaller step size, which makes this time step more expensive.
	Finally, with around $21$ seconds per time step and a larger number of time steps, FT-err is significantly more expensive than FT-dt-u.
	
\section{Conclusions}\label{sec:conclusion}

The two numerical challenges, to derive a feedback control for problems with moving interfaces or free boundaries and to reliably simulate the resulting non-linear closed-loop system,
are tackled successfully in this work.
In detail, we apply the LQR approach and implement non-autonomous versions of the BDF methods and splitting schemes to solve the arising DREs with time-dependent coefficients.
Overall, the splitting methods, in their current state, 
are not competitive for the non-autonomous case.
Other advantages of the splitting methods carry over 
to the non-autonomous case, 
like better predictability of the accuracy 
due to a more uniform error distribution
and less sensitivity to several problem-parameters.
Besides the shorter runtimes, the BDF methods 
can handle more general time-dependencies
in the coefficient matrices.
They are not restricted to the case that, e.g., 
$\cA\time = \alpha\time \bar{\cA}$
and $\cM\time = \mu\time \bar{\cM}$ with constant matrices $\bar{\cA}$ 
and $\bar{\cM}$.
This condition does, e.g., not hold for problems with a time-dependent domain
resulting from moving boundaries or interfaces.
Thus, we exclusively use the BDF methods for the next example, which incorporates a moving interface.

In addition, we have established that the fractional-step-theta scheme with our time-adaptive strategy
can simulate the closed-loop systems with different feedback controls reliably and efficiently.
Namely, our time-adaptivity is tuned to strongly varying controls and prevents blow-ups that can 
occur with other time stepping schemes.
Furthermore, it requires significantly less time steps and runtime compared to classical error-based time-adaptivity.

Forthwith, we observed that the non-autonomous BDF methods and time-adaptive fractional-step-theta scheme, 
in interaction with each other,
are well suited for non-linear large-scale control-problems with moving interfaces.
They perform efficiently and reliably with respect to various problem
parameters and problem dimensions.
After all, despite several stages of approximation and a linearization for the LQR approach, the computed feedback control
stabilizes the non-linear problem successfully.

Besides these results, the runtime performance of the BDF methods, 
as well as the splitting schemes, 
can potentially be improved significantly
by exploring alternative ways to solve the underlying subproblems.
For example, regarding the ARE also available are RADI and projection based
methods like EKSM and RKSM\@.
However, for large problems, a clear advantage of the Newton ADI used in this
work in terms of memory requirements is that only the feedback gain matrices can
be accumulated, thus avoiding the formation of the low-rank solution factors. 
In this case, the restrictions on time-dependent coefficients for splitting schemes could also be removed.

\bibliographystyle{siamplain}
\bibliography{references.bib}

\clearpage
\clearpage
\end{document}


\section{Supplementary Material}	
	\begin{figure}[h]
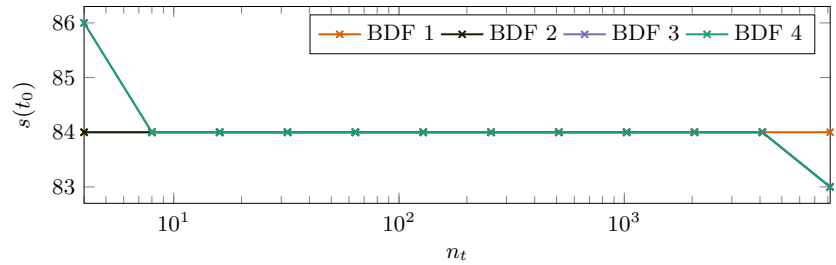

        \vspace{0em}\includetikz{Rail_ranks}
        \caption{\nLtzero{} (numerical rank of the solution) for BDF for different \ntimestep{} (no.\ of time steps), steel profile example.}%
        \label{fig:timeRailranks}
	\end{figure}
	
	\begin{figure}[h]
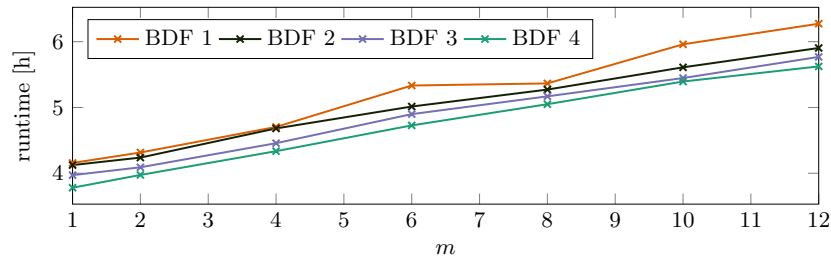

         \vspace{0em}\includetikz{Stefan_times_nb}
         \caption{Runtime of BDF for different \nin{} (columns in \cB\time), Stefan problem example.}%
         \label{fig:timeStefannb}
	\end{figure}

	\begin{figure}[h]
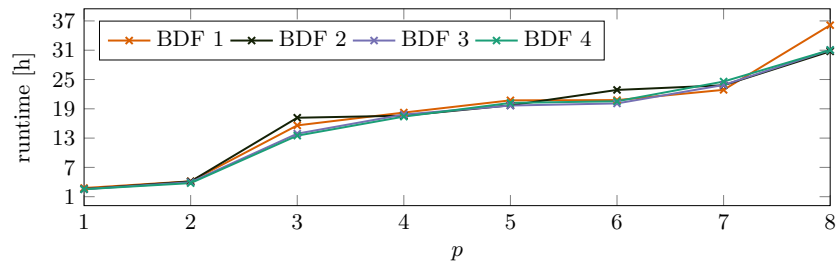

         \vspace{0em}\includetikz{Stefan_times_nc}
         \caption{Runtime of BDF for different \nout{} (rows in \cC\time), Stefan problem example.}%
         \label{fig:timeStefannc}
	\end{figure}

	\begin{figure}[h]
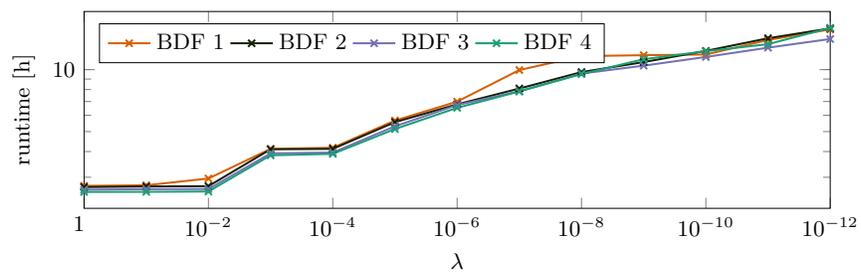

		 \vspace{0em}\includetikz{Stefan_times_control_weight_t_401}
		 \caption{Runtime of BDF for different \weight{} (weight in cost functional), Stefan problem example.}%
		 \label{fig:timeStefancontrolweights}
	\end{figure}

	\begin{figure}[h]
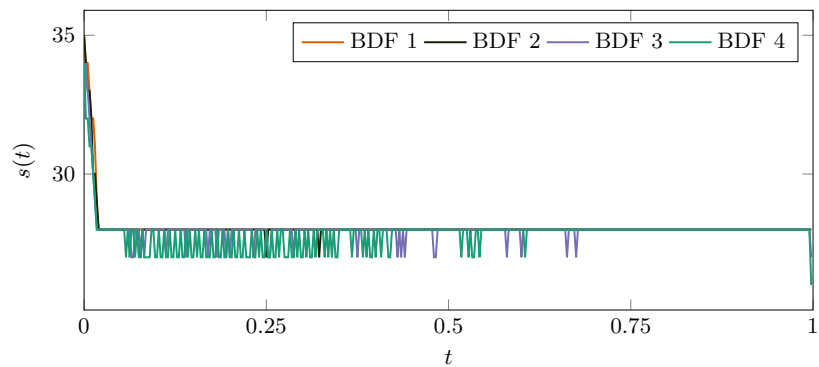

		 \vspace{0em}\includetikz{Stefan_ranks}
		 \caption{Numerical rank of the solution \nLt{} for BDF 1 to 4, Stefan problem example. (used truncation tolerance: machine precision times square of largest singular value)}%
		 \label{fig:Stefanranks}
	\end{figure}